\documentclass[journal]{IEEEtran}
\usepackage{pifont}
\usepackage{}
\usepackage{amsfonts}    
\usepackage{amsmath} 
\usepackage{amssymb}  
\usepackage{mathrsfs}
\usepackage{color,cite}
\usepackage{float} 

\usepackage{subfigure}
\usepackage{graphicx}
\usepackage{epstopdf}
\usepackage{algorithm}
\usepackage{algorithmicx}
\usepackage{algpseudocode}

\usepackage{multicol} 

\begin{document}

\title{Model-free Reinforcement Learning for ${H_{2}/H_{\infty}}$ Control of Stochastic Discrete-time Systems}

\author{
Xiushan Jiang,  ~Li Wang, ~Dongya Zhao,${}^*$   and~Ling Shi

\thanks{ ${}^*$ Corresponding author.}
\thanks{
This work was supported by the National Natural Science Foundation of China (Grant Nos.62103442, 62373229), the Natural Science Foundation of Shandong Province (No.ZR2021QF080), and the Fundamental Research Funds for the Central Universities (Grant No.23CX06024A).
} 
\thanks{
X. Jiang, L. Wang and D. Zhao are with the College of New Energy, China University of Petroleum (East China), Qingdao
266580, China.(e-mail: jiangxsjy@163.com; s21150037@s.upc.edu.cn; dyzhao@upc.edu.cn.)}
\thanks{
Ling Shi is with
the Department of Electronic and Computer Engineering, Hong Kong University of Science and Technology,
Hong Kong, China.(e-mail: eesling@ust.hk.)}

}

%
%

\markboth{Journal of \LaTeX\ Class Files}%
{Shell \MakeLowercase{\textit{et al.}}: Bare Demo of IEEEtran.cls for Journals}
%



\maketitle

\begin{abstract}
This paper proposes a reinforcement learning (RL) algorithm for infinite horizon $\rm {H_{2}/H_{\infty}}$ problem in a class of stochastic discrete-time systems, rather than using a set of coupled generalized algebraic Riccati equations (GAREs). The algorithm is able to  learn the optimal control policy for the system even when its parameters are unknown. Additionally, the paper explores the effect of detection noise as well as the convergence of the algorithm, and shows that the control policy is admissible after a finite number of iterations. The algorithm is also able to handle multi-objective control problems within stochastic fields. Finally, the algorithm is applied to the F-16 aircraft autopilot with multiplicative noise.
\end{abstract}

\newtheorem{theorem}{Theorem}
\newtheorem{lemma}{Lemma}
\newtheorem{corollary}{Corollary}
\newtheorem{proposition}{Proposition}
\newtheorem{definition}{Definition}
\newtheorem{remark}{Remark}
\newtheorem{example}{Example}
\newtheorem{assumption}{Assumption}
\newtheorem{problem}{Problem}

\setlength{\arraycolsep}{2pt}


%
\IEEEpeerreviewmaketitle

\section{Introduction}
{M}{ixed}  $H_{2}/H_{\infty}$ control theory combines both $H_{2}$ control, which optimizes system performance\cite{li2022stochastic}, and $H_{\infty}$ control, which ensures system robustness\cite{huzhipei}. This allows the system to resist external perturbations and minimize a quadratic $H_{2}$ cost under external perturbed inputs. This theory has gained attention from scholars and has been widely used in deterministic systems{ \cite{Doyle1988StatespaceST,khosrowjerdi2004mixed,2002A,LiRobust2018}. In recent years, a body of theoretical research on stochastic differentials and stochastic differences has been gradually constructed and refined \cite{MaoStochastic2007}, and many system control methods have been explored and applied under this new system\cite{li2020mean,li2021stochastic,Yuanpractically2023,DongFlight2020,ZongL12021}. As an advantageous control method, the research on mixed $H_{2}/H_{\infty}$ control has also been extended to the stochastic fields.} For example, Chen and Zhang \cite{chen2004stochastic} solved continuous stochastic $H_{2}/H_{\infty}$ control with state-dependent noises using game theory. They presented a set of related coupled Generalized Algebraic Riccati Equations (GAREs) and investigated the necessary and sufficient conditions for the existence of the solution. In the discrete-time case, Muradorea and Piccib \cite{muradore2005mixed} solved $H_{\infty}$ control for stochastic systems based on the study of deterministic mixed $H_{2}/H_{\infty}$ control, while Zhang et al. \cite{zhang2008infinite,zhang2007stochastic} discussed finite and infinite horizon stochastic $H_{2}/H_{\infty}$ control with $(x, v)$-dependent noises and studied the recursion algorithm for solving the coupled GAREs. Since the GAREs involved are coupled and nonlinear, solving the mixed $H_{2}/H_{\infty}$ problem is a challenging task. Numerical and iterative methods have been widely used to solve the coupled GAREs, but this requires advance knowledge of the system dynamics, which may not hold in practice.

Reinforcement learning (RL) is a machine learning method that uses interactions with the environment to make corrections to actions and allows for real-time optimization of control in a system. It has practical applications {in several fields such as energy, gaming, aviation, and robotics \cite{TongDynamic2023,khan2012reinforcement,KonarADeterministic2013}}. While research on RL in deterministic systems has become more mature over time, several noteworthy contributions have been made. For example, Al-Tamimi et al. \cite{al2007model} proposed a model-free RL method for solving the $ H_{\infty} $ problem, and  Wei et al. \cite{wei2015value} conducted a more in-depth study of the traditional Value Iteration (VI) method, offering insightful conclusions {for the case where} the initial value function is non-zero. Furthermore, Q-learning algorithms have been leveraged to solve a multi-objective nonzero-sum game problem for differential systems \cite{zhang2021q}. To address errors introduced by exploration noise, Kiumarsi \cite{kiumarsi2017h} developed an off-policy algorithm for discrete systems that solves the $ H_{\infty} $ problem.

In recent years, the application of RL algorithms to the optimal control of stochastic systems has gained significant attention {\cite{WuDepth2019,WeiContinuous2023, tancheng1, tancheng2}}. Numerous researchers have examined RL for stochastic systems, particularly for solving the stochastic linear quadratic (SLQ) optimal control problems with RL algorithms. Wang et al. \cite{wang2016infinite,wang2018stochastic} developed an adaptive iterative algorithm that allows for partially model-free solutions of SLQ optimal control. Pang and Jiang \cite{pang2022reinforcement} improved the traditional policy iterative algorithm for continuous systems based on the least squares method. Li et al. \cite{li2022stochastic} solved the optimal SLQ problem for a continuous time system with $(x, u)$-dependent noise using  {RL algorithms}. These algorithms have been shown to effectively solve the optimal control problem in stochastic systems even with partially unknown models. Additionally, recent research has focused on completely model-free algorithms. For example, Lai et al. \cite{lai2023model} developed online and offline Q-learning algorithms to directly deal with discrete stochastic systems with both multiplicative and additive noise, demonstrating their superior performance through simulation experiments. However, existing research on RL in stochastic systems only focuses on single-objective problems, which cannot solve complex coupled GAREs for multi-objective cases. Although no research has been conducted on applying {RL} algorithms to stochastic $H_{2}/H_{\infty}$ problems, the results from the aforementioned literature serve as the basis for this paper.

This paper aims to develop a model-free RL algorithm to solve the mixed $H_{2}/H_{\infty}$ problem for discrete-time stochastic systems. The contributions of this paper are as follows.
\begin{itemize}
	\item[1.]  We adopt RL algorithms to solve the stochastic mixed $ H_{2}/H_{\infty} $ control problem for the first time, instead of numerical methods \cite{chen2004stochastic,zhang2008infinite,zhang2007stochastic}, which is significant for the case where the system dynamics model is inaccessible in practice.
	\item[2.]  We explore the issue about the admissibility of control strategies, which is a complement to the research on Q-learning algorithms in \cite{wang2018stochastic}.
	\item[3.]  Unlike the single-objective optimal control problem in { \cite{wang2016infinite,wang2018stochastic,li2022stochastic,lai2023model,pang2022reinforcement}}, we consider two optimization objectives to ensure both robustness and optimality, which is much more challenging to tackle. Compared to \cite{zhang2021q}, we extend Q-learning under multiple objectives to stochastic systems.
\end{itemize}

The remainder of this paper is organized as follows. In Section \uppercase\expandafter{\romannumeral2}, some preliminaries are given, and the $ H_{2}/H_{\infty} $ problem for stochastic discrete-time systems is described. In Section \uppercase\expandafter{\romannumeral3}, the Q-learning algorithm is developed and the convergence is proved. A simulation example is given in Section \uppercase\expandafter{\romannumeral4} to demonstrate the effectiveness of the proposed algorithm. Finally, the conclusions are given in Section \uppercase\expandafter{\romannumeral5}.

{\em Notations}: $ M>0 $  $ (M \geq 0) $ means that $ M $ is a real symmetric positive definite (semi-positive definite) matrix. $ M^{T} $ is the transpose of the $ M $. $ Tr(M) $ is the trace of the square matrix $ M $. $ \mathcal{R}^{m}$ represents a space of $ m $-dimensional real vectors, and $ \mathcal{R}^{m\times n} $ denotes a set of $ m\times n $ real matrices. $ E(x) $ is used to express the mathematical expectation of the random variable $ x $. $ \left\|\cdot \right\| $ is the Euclidean norm of $ n $-dimensional vectors; $ M \otimes N $ is the Kronecker product of $ M $ and $ N $. $ \mathcal{N}_{k}:=\left\lbrace 0,1,2,\cdots,k\right\rbrace  $ and $ \mathcal{N}_{+}:=\left\lbrace 0,1,2,3,\cdots\right\rbrace  $.
\section{Problem formulation and preliminary}
In this paper, we consider the following stochastic discrete linear time-invariant (SDLTI) system with state- and disturbance-dependent noise
\begin{equation}\label{system}
	\left\lbrace
	\begin{split}
		&x_{k+1}= A_{1}x_{k}+B_{1}u_{k}+C_{1}v_{k} +\left( A_{2}x_{k}+C_{2}v_{k}\right) \omega_{k},\\
		&{y_{k}}=\left[ \begin{matrix}Cx_{k}\\Du_{k}\end{matrix} \right] , D^{T}D=I,
	\end{split}\right.
\end{equation}
where $ x_{k}\in\mathcal{R}^n $ denotes the system state with the any deterministic initial value $ x_{0} $. The time index is $ k\in \mathcal{N}_{+} $. $ u_{k}\in\mathcal{R}^{m_{1}} $, $ v_{k}\in\mathcal{R}^{m_{2}} $ and $ y_{k}\in\mathcal{R}^{m} $ are the control input, external disturbance input and controlled output, respectively. $ A_{1}, A_{2}\in\mathcal{R}^{n\times n} $, $ B_{1}\in\mathcal{R}^{n\times m_{1}} $ and $ C_{1}, C_{2}\in\mathcal{R}^{n\times m_{2}} $ are time-invariant system dynamical matrices. $ \left\lbrace \omega_{k}, k \in \mathcal{N}_{+} \right\rbrace $ represents an independent white noise sequence defined on a complete filtered probability space ($\Omega, \mathcal{F}, \mathcal{P} $), which is a wide sense stationary process satisfying $ E\left( \omega_{t}\right) = 0 $ and $ E\left( \omega_{s}\omega_{k}\right)  = \delta_{sk} $. $ \delta_{sk} $ denotes the Kronecker delta so that
\begin{equation}
	\left\lbrace
	\begin{split}
		&E\left( \omega_{s}\omega_{k}\right)=1~~~~~{\rm if}~~s=k,\\
		&E\left( \omega_{s}\omega_{k}\right)=0~~~~~{\rm if}~~s\neq k.
	\end{split}\right. \nonumber
\end{equation}
$ \mathcal{F}_{k} := \sigma\left( \omega_{s} : s\in \mathcal{N}_{k}\right)  $ is the $ \sigma $-algebra generated by $ \omega_{s}, s \in \mathcal{N}_{+} $. Let $ l_{\omega }^{2}(\mathcal{N}_{+}, \mathcal{R}^{n}) $ be a space consisting of all non-anticipative stochastic processes $ \phi_{k}\in \left\lbrace \mathcal{N}_{+} \times \Omega \rightarrow \mathcal{R}^{n} \right\rbrace $, such that $ \phi_{k} $ is $ \mathcal{F}_{k-1} $ measurable, in which we define $ \mathcal{F}_{-1} = \left\lbrace \emptyset, \Omega \right\rbrace  $  for a constant initial value $ \phi_{0} $. The $ l^{2} $-norm of $ \phi_{k}\in l_{\omega}^{2}(\mathcal{N}_{+}, \mathcal{R}^{n}) $ satisfies $$ \left\| \phi\right\|_{l_{\omega}^{2}(\mathcal{N}_{+}, \mathcal{R}^{n})} :=\left( \sum_{k=0}^{\infty}E\left\| \phi_{k}\right\|^{2} \right)^{\frac{1}{2}} < +\infty. $$

Before describing the problem to be solved, some necessary definitions and lemmas are given.
\begin{definition}
	\cite{zhang2017stochastic}~The following SDLTI system
	$$
	x_{k+1}=A_{1}x_{k}+A_{2}x_{k} { \omega_{k}}
	$$
	is said to be asymptotically stable in the mean square (ASMS) sense if the corresponding system state satisfies $$\lim_{k\to\infty}E\left(x_{k}^{T}x_{k} \right)=0 $$ for any initial state $ x_{0}\in \mathcal{R}^{n} $.
\end{definition}
\begin{definition}
	\cite{zhang2017stochastic}~The control policy $ u_{k}=K_{2}x_{k}  $ is said to be admissible if the following closed-loop system
	$$ x_{k+1}=\left( A_{1}+B_{1}K_{2}\right) x_{k}+A_{2}x_{k}{ \omega_{k}} $$
	is ASMS.
\end{definition}
\begin{definition}
	\cite{zhang2017stochastic}~The following system
	\begin{equation}
		\left\lbrace \begin{split}
			&x_{k+1}=A_{1}x_{k}+A_{2}x_{k}{\omega_{k}},\\
			&{y_{k}}=Cx_{k}, x_{0} \in \mathcal{R}^{n}, k \in \mathcal{N}_{+}
		\end{split}
		\right. \nonumber
	\end{equation}
	is exactly observable if {$ y_{k} \equiv 0, a.s., \forall k \in  \mathcal{N}_{+} \Rightarrow x_{0}=0 . $ }	
\end{definition}
\begin{definition}
	The perturbation operator of system (\ref{system}) is defined as $ \mathcal{L}(v) := \left[ \begin{matrix}C{x_{k}}\\Du^{*}\end{matrix} \right] $ with its norm
	\begin{equation}
		\begin{split}
			\left\| \mathcal{L}\right\| :=\sup\limits_{v\in l_{\omega}^{2}(\mathcal{N}_{+},\mathcal{R}^{m_{2}}),\atop v \neq 0,x_{0}=0}\dfrac{\left(  \sum\limits_{k=0}^{\infty}E(\left\| Cx_{k}\right\|^{2} + \left\| {u_{k}^{*}}\right\|^{2} )\right)  ^{\frac{1}{2}} }{\left(  \sum\limits_{k=0}^{\infty}Ev_{k}^{T}v_{k} \right)  ^{\frac{1}{2}}} ,
		\end{split}
		\nonumber
	\end{equation}
	{where $x_{k}$ is the system state generated by the optimal control strategy $u_{k}^{*}\in l_{\omega}^{2}(\mathcal{N}_{+},\mathcal{R}^{m_{1}})$ with $x_{0}=0$.}
\end{definition}
\begin{lemma}\label{LemLyapunov}
	\cite{jiang2019stability}~For the SDLTI system $$ x_{k+1}=A_{1}x_{k}+A_{2}x_{k}{ \omega_{k}},$$ it is ASMS if and only if the Lyapunov-type inequality $$ A_{1}^{T}PA_{1}+A_{2}^{T}PA_{2}-P < 0 $$ has at least one solution $ P > 0 $.
\end{lemma}

To facilitate the subsequent derivation in this paper, we make the following assumptions without loss of generality.
\begin{assumption}
	For the system (\ref{system}), there exists at least one admissible control.
\end{assumption}	
\begin{assumption}
	The system (\ref{system}) is exactly observable.
\end{assumption}
\begin{remark}
	Exact observability can be verified and guaranteed using different criteria \cite{zhang2017stochastic}. In most cases, exact observability can be directly inferred from $C^TC>0$.
\end{remark}

We consider the system \eqref{system} with a prescribed disturbance attenuation level $ \gamma>0 $. Define two associated performance indices as
\begin{equation}
	\left\lbrace \begin{split}
		J_{1}\left( x_{0},u,v\right) &= \sum\limits_{k=0}^{\infty} r_{1}\left(  x_{k},u_{k},v_{k}\right) \\&= \sum\limits_{k=0}^{\infty} E\left[ \gamma^{2} \left\|  v_{k}\right\|^{2} -\left\|  y_{k}\right\| ^{2} \right] ,\\
		J_{2}\left( x_{0},u,v\right) &= \sum\limits_{k=0}^{\infty} r_{2}\left(  x_{k},u_{k},v_{k}\right)= \sum\limits_{k=0}^{\infty} E\left\|  y_{k}\right\| ^{2}.
	\end{split}\right.
\end{equation}
We set the given performance function to be associated with the controlled output matrix with $ Q= C^{T}C \in \mathcal{R}^{n\times n}$, whose value is able to be obtained in advance. And
\begin{equation}
	\begin{split}
		r_{1}\left(  x_{k},u_{k},v_{k}\right)=&E\left(\gamma^{2} v_{k}^{T}v_{k}-x_{k}^{T}Qx_{k}-u_{k}^{T}u_{k} \right),  \\
		r_{2}\left(  x_{k},u_{k},v_{k}\right)=&E\left(x_{k}^{T}Qx_{k}+u_{k}^{T}u_{k} \right).
	\end{split}
	\nonumber
\end{equation}
\begin{definition}
	$v^{*}_{k}$ is said to be the worst disturbance if $v^{*}_{k}$  minimizes the performance function $ J_{1}(x_{k},u_{k},v_{k}) $ for all disturbances $ v_{k} $.
\end{definition}

Based on the previous system description and performance functions, we formally present the problem to be solved in this paper, which is named the infinite horizon stochastic $ H_{2}/H_{\infty} $ problem.
\begin{problem}\label{problem}
	(Infinite horizon stochastic $H_{2}/H_{\infty}$ problem). For system (\ref{system}), the  infinite horizon stochastic $ H_{2}/H_{\infty} $ problem is to find the optimal feedback control strategy $ u_{k}^{*}\in l_{\omega}^{2}(\mathcal{N}_{+}, \mathcal{R}^{m_{1}}) $ and the worst-case disturbance $ v_{k}^{*} \in l_{\omega}^{2}(\mathcal{N}_{+}, \mathcal{R}^{m_{2}}) $ such that
	\begin{itemize}
		\item[1)]  $u_{k}^{*}$ makes the system (\ref{system}) ASMS in the absence of $v_{k}${.}
		\item[2)]  The inequality $\left\| \mathcal{L}\right\|<\gamma$ always holds for any disturbances $ v_{k} $.
		\item[3)]  $u_{k}^{*}$ minimizes the output energy $ J_{2}(x_{k},u_{k},v^{*}_{k})$ when considering the worst-case disturbance $v^{*}_{k}$.
	\end{itemize}
\end{problem}

Since the system (\ref{system}) is exactly observable, the coupled matrix-valued equations \eqref{P1Riccati}-\eqref{P2Riccati} given by Theorem 3.8 in \cite{zhang2017stochastic} have a solution $ (P_{1}, K_{1}; P_{2}, K_{2}) $ with $  P_{1} < 0 $ and $ P_{2} > 0 $  so that the previous Problem \ref{problem} is solvable \cite{zhang2008infinite}.
\begin{align}
	&-P_{1}+\left(A_{1}+B_{1}K_{2} \right)^{T}P_{1}\left(A_{1}+B_{1}K_{2} \right)-Q \nonumber\\
	& +A_{2}^{T}P_{1}A_{2}-K_{2}^{T}K_{2}-M_{1}\Delta_{1}^{-1}M_{1}^{T}=0,\label{P1Riccati} \\
	&-P_{2}+\left(A_{1}+C_{1}K_{1} \right)^{T}P_{2}\left(A_{1}+C_{1}K_{1} \right) +Q\nonumber\\
	&+\left(A_{2}+C_{2}K_{1}\right)^{T}P_{2}\left( A_{2}+C_{2}K_{1}\right)
	-M_{2}\Delta_{2}^{-1}M_{2}^{T}=0,\label{P2Riccati}
\end{align}
Then the solution to Problem \ref{problem} associated with \eqref{P1Riccati}-\eqref{P2Riccati} is formulated as $ u^{*}_{k}={ K_{2}}x_{k}$ and $v^{*}_{k}=K_{1}x_{k} $, where
\begin{align}
	K_{1}=-\Delta_{1}^{-1}M_{1}^{T} ,\label{K1Riccati}\\
	K_{2}=-\Delta_{2}^{-1}M_{2}^{T}.\label{K2Riccati}
\end{align}
In the above equations (\ref{P1Riccati})-(\ref{K2Riccati}), the matrices $ M_{1} ,\Delta_{1},M_{2},\Delta_{2}$ are given as
\begin{equation}
	\begin{split}
		& M_{1}=\left(A_{1}+B_{1}K_{2} \right)^{T}P_{1}C_{1}+A_{2}^{T}P_{1}C_{2} ,\\
		&\Delta_{1}=\gamma^{2}I+C_{2}^{T}P_{1}C_{2}+C_{1}^{T}P_{1}C_{1}>0,\\
		&M_{2}=\left(A_{1}+C_{1}K_{1} \right)^{T}P_{2}B_{1},\\
		& \Delta_{2}=I+B_{1}^{T}P_{2}B_{1}>0,
	\end{split}
	\nonumber
\end{equation}
and $ rank\left( \begin{bmatrix}\Delta_{1}&C_{1}^{T}P_{1}B_{1}\\B_{1}^{T}P_{2}C_{1}&\Delta_{2}		\end{bmatrix}\right) =m_{1}+m_{2} $. The optimal performance corresponding to the system \eqref{system} can be expressed as $ J_{1}(x_{0},u^{*},v^{*})=x_{0}^{T}P_{1}x_{0}\leq J_{1}(x_{0},u^{*},v) $ and $ J_{2}(x_{0},u^{*},v^{*})=x_{0}^{T}P_{2}x_{0}\leq J_{2}(x_{0},u,v^{*})$.

Based on {the} above discussion, combined with (\ref{K1Riccati}) and (\ref{K2Riccati}), $ K_{1}$ and $K_{2} $ can be converted as
$$  \begin{bmatrix}\Delta_{1}&C_{1}^{T}P_{1}B_{1}\\B_{1}^{T}P_{2}C_{1}&\Delta_{2}		\end{bmatrix}\begin{bmatrix}	K_{1}\\K_{2}\end{bmatrix}
+\begin{bmatrix}C_{1}^{T}P_{1}A_{1}+ C_{2}^{T}P_{1}A_{2}\\B_{1}^{T}P_{2}A_{1}\end{bmatrix}=0. $$
With invertibility of the first matrix in the above equation, it follows that
\begin{align} \label{K1K2Riccati}
	\begin{bmatrix}	K_{1}\\K_{2}\end{bmatrix}=&-{\begin{bmatrix}
			\gamma^{2}I+C_{2}^{T}P_{1}C_{2}+C_{1}^{T}P_{1}C_{1}&C_{1}^{T}P_{1}B_{1}\\B_{1}^{T}P_{2}C_{1}&I+B_{1}^{T}P_{2}B_{1}	
		\end{bmatrix}^{-1}} \nonumber \\
	&\times \begin{bmatrix}	C_{1}^{T}P_{1}A_{1}+ C_{2}^{T}P_{1}A_{2}\\B_{1}^{T}P_{2}A_{1}\end{bmatrix}.
\end{align}

Now we obtain the solution to the Problem \ref{problem} of system (\ref{system}). $  u_{k}^{*} = K_{2}x_{k} $ is the desired $ H_{\infty} $ control, while $  v_{k}^{*}=K_{1}x_{k} $ is the corresponding worst-case disturbance. It is noteworthy that the solution to Problem \ref{problem} is related to the coupled GAREs \eqref{P1Riccati}-\eqref{P2Riccati} and depends on the system dynamical matrices. On the one hand, it is challenging to solve the coupled matrix-valued equations directly. On the other hand, the solution requires prior knowledge of the system dynamic parameters, which places higher demands on system \eqref{system}. To avoid these problems, we will develop a model-free algorithm for Problem \ref{problem} in next section.
\section{Model-free RL for stochastic $ H_{2}/H_{\infty} $ problem}
Considering the difficulties in obtaining the solution to stochastic $ H_{2}/H_{\infty} $ problem, we will construct a model-free RL algorithm to solve the coupled GAREs without prior knowledge of the system dynamics parameters in this section.

\subsection{Q-function Setup for Stochastic Systems}
We develop the model-free RL algorithm based on Q-functions. Setting the initial time to be $ k $, the cost function can be written as
\begin{align}
	&J_{1}(x_{k},u_{k},v_{k}) = r_{1}\left(  x_{k},u_{k},v_{k}\right)+\sum\limits_{t=k+1}^{\infty}r_{1}\left(  x_{t},u_{t},v_{t}\right) \nonumber\\
	&=r_{1}\left(  x_{k},u_{k},v_{k}\right)+J_{1}(x_{k+1},u_{k+1},v_{k+1}),\label{J1:k}\\
	&J_{2}(x_{k},u_{k},v_{k})=r_{2}\left(  x_{k},u_{k},v_{k}\right)+\sum\limits_{t=k+1}^{\infty}r_{2}\left(  x_{t},u_{t},v_{t}\right) \nonumber\\
	&=r_{2}\left(  x_{k},u_{k},v_{k}\right)+J_{2}(x_{k+1},u_{k+1},v_{k+1}).\label{J2:k}
\end{align}
Substituting the optimal solution $ (u^{*}_{k}, v^{*}_{k}) $ into equations (\ref{J1:k}) and (\ref{J2:k}), the corresponding optimal performance of the system (\ref{system}) can be found as follows
\begin{align}
	J_{1}^{*}(x_{k}):=& J_{1}(x_{k},u_{k}^{*},v_{k}^{*})\nonumber\\
	=&\min \limits_{v_{k}} J_{1}(x_{k},u_{k}^{*},v_{k})=E\left( x_{k}^{T}P_{1}x_{k}\right) ,\\
	J_{2}^{*}(x_{k}):=& J_{2}(x_{k},u_{k}^{*},v_{k}^{*})\nonumber\\
	=&\min \limits_{{u_{k}}} J_{2}(x_{k},u_{k},v_{k}^{*})=E\left( x_{k}^{T}P_{2}x_{k}\right) ,
\end{align}
where $ P_{1} $ and $ P_{2} $ are the solutions of equations (\ref{P1Riccati}) and (\ref{P2Riccati}), respectively.

It is desired to find the optimal solution $ (u_{k}^{*}, v_{k}^{*}) $. According to Bellman's principle of optimality, the optimal cost function should satisfy the Hamilton-Jacobi-Bellman (HJB) equations
\begin{align}
	J_{1}^{*}(x_{k}) = \min \limits_{v_{k}} \left[ r_{1}\left(  x_{k},u_{k},v_{k}\right)+J_{1}^{*}(x_{k+1})\right],\label{J1*}\\
	J_{2}^{*}(x_{k}) = \min \limits_{u_{k}}\left[ r_{2}\left(  x_{k},u_{k},v_{k}\right)+J_{2}^{*}(x_{k+1})\right],\label{J2*}
\end{align}
where $ u_{k} = u_{k}^{*} $ is in the equation (\ref{J1*}), and $ v_{k} = v_{k}^{*} $ is in the equation (\ref{J2*}).

In order to develop an online model-free algorithm, we first give the definition of Q-functions $ Q_{1}^{*}$ and $Q_{2}^{*} $ according to (\ref{J1*}) and (\ref{J2*}). Define the vector $ z_{k}$ as $ z_{k}:= \begin{bmatrix}x_{k}^{T}&u_{k}^{T}&v_{k}^{T}\end{bmatrix}^{T}  $. Then
\begin{align}
	&Q_{1}^{*}(x_{k},u_{k},v_{k})\nonumber
	:= r_{1}\left(  x_{k},u_{k},v_{k}\right)+J_{1}^{*}(x_{k+1})\nonumber \\
	=&E\left(\gamma^{2} v_{k}^{T}v_{k}-x_{k}^{T}Qx_{k}-u_{k}^{T}u_{k} \right)+E\left( x_{k+1}^{T}P_{1}x_{k+1}\right) \nonumber \\
	=&E \begin{bmatrix}	 x_{k}\\u_{k}\\v_{k} \end{bmatrix}^{T}
	\left( \begin{bmatrix}-Q&0&0\\0&-I&0\\0&0&\gamma^{2}I\end{bmatrix}+
	\begin{bmatrix}
		A_{1}^{T}\\B_{1}^{T}\\C_{1}^{T}
	\end{bmatrix}
	P_{1}\begin{bmatrix}
		A_{1}^{T}\\B_{1}^{T}\\C_{1}^{T}
	\end{bmatrix}^{T}\right. \nonumber \\
	&\left. +\begin{bmatrix}
		A_{2}^{T}\\0\\C_{2}^{T}
	\end{bmatrix}
	P_{1}		\begin{bmatrix}
		A_{2}^{T}\\0\\C_{2}^{T}
	\end{bmatrix}^{T} \right)
	\begin{bmatrix}	 x_{k}\\u_{k}\\v_{k} \end{bmatrix} \nonumber \\
	=&E \left( z_{k}^{T}H_{1}z_{k}\right).
\end{align}
For convenience, we set $ p := n+m_{1}+m_{2}  $ and $ p \in N_{k}  $. The matrix $ H_{1}\in \mathcal{R}^{p\times p} $ is symmetric, which can be rewritten as
\begin{equation}\label{H1}
	\begin{split}
		H_{1}=\begin{bmatrix}
			\Gamma_{1}(1,1) & \Gamma_{1}(1,2) & \Gamma_{1}(1,3)\\
			\Gamma^{T}_{1}(1,2) & \Gamma_{1}(2,2) & \Gamma_{1}(2,3)\\
			\Gamma^{T}_{1}(1,3) & \Gamma^{T}_{1}(2,3) & \Gamma_{1}(3,3)
		\end{bmatrix}
	\end{split}
\end{equation}
with
\begin{equation}
	\begin{split}
		&\Gamma_{1}(1,1)=A_{1}^{T}P_{1}A_{1}+A_{2}^{T}P_{1}A_{2}-Q,\\
		&\Gamma_{1}(1,2)=A_{1}^{T}P_{1}B_{1},\\
		&\Gamma_{1}(1,3)=A_{1}^{T}P_{1}C_{1}+A_{2}^{T}P_{1}C_{2},\\
		&\Gamma_{1}(2,2)=B_{1}^{T}P_{1}B_{1}-I,\\
		&\Gamma_{1}(2,3)=B_{1}^{T}P_{1}C_{1},\\
		&\Gamma_{1}(3,3)=C_{1}^{T}P_{1}C_{1}+C_{2}^{T}P_{1}C_{2}+\gamma^{2}I.
	\end{split}
	\nonumber
\end{equation}
Similarly, it can be seen that
\begin{equation}
	\begin{split}
		Q_{2}^{*}(x_{k},u_{k},v_{k}):= &r_{2}\left(  x_{k},u_{k},v_{k}\right)+J_{2}^{*}(x_{k+1})\\
		=&E \left( z_{k}^{T}  H_{2} z_{k}\right) ,
	\end{split}
\end{equation}
where the symmetric matrix $ H_{2}\in \mathcal{R}^{p\times p} $ can be expressed as
\begin{equation}\label{H2}
	\begin{split}
		&H_{2}=\begin{bmatrix}
			\Gamma_{2}(1,1) & \Gamma_{2}(1,2) & \Gamma_{2}(1,3)\\
			\Gamma^{T}_{2}(1,2) & \Gamma_{2}(2,2) & \Gamma_{2}(2,3)\\
			\Gamma^{T}_{2}(1,3) & \Gamma^{T}_{2}(2,3) & \Gamma_{2}(3,3)
		\end{bmatrix}
	\end{split}
\end{equation}
with
\begin{equation}
	\begin{split}
		&\Gamma_{2}(1,1)=A_{1}^{T}P_{2}A_{1}+A_{2}^{T}P_{2}A_{2}+Q,\\
		&\Gamma_{2}(1,2)=A_{1}^{T}P_{2}B_{1},\\
		&\Gamma_{2}(1,3)=A_{1}^{T}P_{2}C_{1}+A_{2}^{T}P_{2}C_{2},\\
		&\Gamma_{2}(2,2)=B_{1}^{T}P_{2}B_{1}+I,\\
		&\Gamma_{2}(2,3)=B_{1}^{T}P_{2}C_{1},\\
		&\Gamma_{2}(3,3)=C_{1}^{T}P_{2}C_{1}+C_{2}^{T}P_{2}C_{2}.
	\end{split}
	\nonumber
\end{equation}
It can be noted that the value function $ Q_{1}^{*}(x_{k},u_{k}, { v_{k}})$, $ Q_{2}^{*}(x_{k},u_{k},{ v_{k}})$ are equal to the optimal value of the cost function $ J_{1}^{*}(x_{k})$,  $ J_{2}^{*}(x_{k})$, respectively, when the control strategy $ u_{k} $ is optimal and the external disturbance $ v_{k} $ is the worst. Then one has
\begin{align}
	J_{1}^{*}(x_{k}) = \min\limits_{v_{k}} Q_{1}^{*}(x_{k},u^{*}_{k},v_{k})=Q_{1}^{*}(x_{k},u^{*}_{k},v^{*}_{k}),\\
	J_{2}^{*}(x_{k}) = \min\limits_{u_{k}} Q_{2}^{*}(x_{k},u_{k},v_{k}^{*})=Q_{2}^{*}(x_{k},u^{*}_{k},v^{*}_{k}).
\end{align}
Thus, we obtain the form of matrix pair $ \left( H_{1}, H_{2}\right)  $ in $ Q_{1}^{*}(x_{k},u_{k},v_{k}) $ and $ Q_{2}^{*}(x_{k},u_{k},v_{k}) $ corresponding to the HJB equations (\ref{J1*}) and (\ref{J2*}) as
\begin{align}
	Q_{1}^{*}(x_{k},u_{k},v_{k}) = r_{1}\left(  x_{k},u_{k},v_{k}\right)+Q_{1}^{*}(x_{k+1},u_{k+1}^{*},v_{k+1}^{*}),\\
	Q_{2}^{*}(x_{k},u_{k},v_{k}) = r_{2}\left(  x_{k},u_{k},v_{k}\right)+Q_{2}^{*}(x_{k+1},u_{k+1}^{*},v_{k+1}^{*}).
\end{align}
Using the internal information of matrices $ H_{1} $ and $ H_{2} $ in (\ref{H1}) and (\ref{H2}), we are ready to re-represent the gains of the $ H_{2}/H_{\infty} $ control policy and disturbance in (\ref{K1K2Riccati}), such as
\begin{equation}\label{K1K2 by H1H2}
	\begin{split}
		\begin{bmatrix}	K_{1}\\K_{2}\end{bmatrix}=-\begin{bmatrix}
			\Gamma_{1}(3,3)&\Gamma_{1}^{T}(2,3)\\\Gamma_{2}(2,3)&\Gamma_{2}(2,2)		
		\end{bmatrix}^{-1} \begin{bmatrix}	\Gamma_{1}^{T}(1,3)\\\Gamma_{2}^{T}(1,2)\end{bmatrix}.
	\end{split}
\end{equation}
From (\ref{K1K2 by H1H2}), it can be seen that the infinite mixed $ H_{2}/H_{\infty} $ solution $ (u_{k}^{*}, v_{k}^{*}) $ is completely determined by the matrix pair $ (H_{1}, H_{2}) $. It is obvious that the matrices $ H_{1} $ and $ H_{2} $ contain information about the dynamics of the system. In the case where $ (H_{1}, H_{2}) $ is known, solving the optimal control policy $ u_{k}^{*} $ and the worst-case perturbation $ v_{k}^{*}$ does not require a system dynamic model. Therefore, the equation (\ref{K1K2 by H1H2}) is the key for the Q-learning algorithm to implement the model-free solution to Problem \ref{problem}.

Next, we present the Q-learning iterative algorithm to obtain the values of matrix pair $ (H_{1}, H_{2}) $ online. The developed algorithm is able to solve the coupled GAREs by collecting trajectories of the stochastic system without knowledge of the dynamical model of system (\ref{system}).
\subsection{Online Implementation of the Model-free Algorithm}
Based on the constructed Q-functions, now we present the derivation and implementation of the Q-learning algorithm. In contrast to the VI algorithm, which can be implemented partially model-free, the Q-learning algorithm is a method for solving optimal control problems without any knowledge of the system dynamics.

In the proposed Q-learning algorithm, we start with initial Q-functions $ Q_{1}^{(0)}(\cdot) $ and $ Q_{2}^{(0)}(\cdot) $ (for convenience, they are set to 0) and $ (u_{k}^{(0)},v_{k}^{(0)}) $. The Q-functions at this point are not necessarily optimal. Then we obtain $ Q_{1}^{(1)}(x_{k},{u_{k}^{(0)},v_{k}^{(0)}})$, $ Q_{2}^{(1)}(x_{k},{u_{k}^{(0)},v_{k}^{(0)}}) $ by setting $ i = 0 $ in (\ref{Q1Q2:i}).
\begin{equation}\label{Q1Q2:i}
	\left\lbrace \begin{split}
		&Q_{1}^{(i+1)}(x_{k},{u_{k}^{(i)},v_{k}^{(i)}}) \\
		=&r_{1}\left(  x_{k},u_{k}^{(i)},v_{k}^{(i)}\right)+E\left( x_{k+1}^{T}P_{1}^{(i)}x_{k+1}\right) \\
		=&r_{1}\left(  x_{k},u_{k}^{(i)},v_{k}^{(i)}\right)+E\begin{bmatrix}	 x_{k+1}\\u_{k+1}^{(i)}\\v_{k+1}^{(i)} \end{bmatrix}^{T}
		H_{1}^{(i)}
		\begin{bmatrix}	 x_{k+1}\\u_{k+1}^{(i)}\\v_{k+1}^{(i)} \end{bmatrix},\\
		&Q_{2}^{(i+1)}(x_{k},{u_{k}^{(i)},v_{k}^{(i)}}) \\
		=&r_{2}\left(  x_{k},u_{k}^{(i)},v_{k}^{(i)}\right)+E\left( x_{k+1}^{T}P_{2}^{(i)}x_{k+1}\right) \\
		=&r_{2}\left(  x_{k},u_{k}^{(i)},v_{k}^{(i)}\right)+E\begin{bmatrix}	 x_{k+1}\\u_{k+1}^{(i)}\\v_{k+1}^{(i)} \end{bmatrix}^{T}
		H_{2}^{(i)}
		\begin{bmatrix}	 x_{k+1}\\u_{k+1}^{(i)}\\v_{k+1}^{(i)} \end{bmatrix}.
	\end{split}\right.
\end{equation}
The generalized greedy algorithm is used in the calculation of Q-function, i.e., the optimization of $ Q^{(i+1)} $ is achieved by adopting (\ref{vu:i}) at the $ i^{th} $ iteration. Thus, $ (u_{k}^{(i)},v_{k}^{(i)}) $ will satisfy
\begin{equation}\label{vu:i}
	\left\lbrace \begin{split}	
		v_{k}^{(i+1)}= \arg \min \limits_{v_{k}} & Q_{1}^{(i+1)}(x_{k},u^{(i+1)}_{k},v_{k})\\
		=\arg \min \limits_{v_{k}} &\left[ r_{1}\left(  x_{k},u_{k}^{(i+1)},v_{k} \right)\right.\\
		& \left. +Q_{1}^{(i)}(x_{k+1},u_{k+1}^{(i+1)},v_{k+1} )\right] ,\\
		u_{k}^{(i+1)}= \arg \min \limits_{u_{k}} & Q_{2}^{(i+1)}(x_{k},u_{k},v_{k}^{(i+1)})\\
		=\arg \min \limits_{u_{k}} &\left[ r_{2}\left(  x_{k},u_{k},v_{k}^{(i+1)}\right) \right.\\
		& \left.+Q_{2}^{(i)}(x_{k+1},u_{k+1},v_{k+1}^{(i+1)})\right] .
	\end{split} \right.
\end{equation}
Based on (\ref{K1K2 by H1H2}), the updated gains $ (K_{1}^{(i)},K_{2}^{(i)}) $ can be expressed by $ H_{1} $ and $ H_{2} $ as
\begin{equation}\label{K1K2:i by H1H2}
	\begin{split}
		\begin{bmatrix}	K_{1}^{(i)}\\K_{2}^{(i)}\end{bmatrix}=-\begin{bmatrix}
			\Gamma_{1}^{(i)}(3,3)&\Gamma_{1}^{(i)T}(2,3)\\\Gamma_{2}^{(i)}(2,3)&\Gamma_{2}^{(i)}(2,2)		
		\end{bmatrix}^{-1} \begin{bmatrix}	\Gamma_{1}^{(i)T}(1,3)\\\Gamma_{2}^{(i)T}(1,2)\end{bmatrix}
	\end{split},
\end{equation}
$$ u_{k}^{(i)}={K_{2}^{(i)}}x_{k}, v_{k}^{(i)}={K_{1}^{(i)}}x_{k}. $$
Then, it follows that
\begin{equation}\label{Pi:Hi}
	\left\lbrace \begin{split}	
		&P_{1}^{(i)}= \begin{bmatrix}I&K_{2}^{(i)T}&K_{1}^{(i)T}	\end{bmatrix}  H_{1}^{(i)} \begin{bmatrix}I&K_{2}^{(i)T}&K_{1}^{(i)T}	\end{bmatrix}^{T},\\
		&P_{2}^{(i)}= \begin{bmatrix}I&K_{2}^{(i)T}&K_{1}^{(i)T}	\end{bmatrix}  H_{2}^{(i)} \begin{bmatrix}I&K_{2}^{(i)T}&K_{1}^{(i)T}	\end{bmatrix}^{T}.
	\end{split}\right.
\end{equation}
\begin{remark}
	It is essential to note that $ (u_{k}^{(i)}, v_{k}^{(i)}) $ is updated synchronously during the iterative process. The update of the action networks is determined only by the information of the matrix pair $ (H_{1}, H_{2}) $ and does not require the model parameters of the system (\ref{system}). According to (\ref{Q1Q2:i}), it is easy to obtain the matrices $ H_{1}$ and $H_{2} $ that depend only on the input information and the system state that are available online. It can be observed that the Q-learning algorithm is a model-free time-forward method for solving the optimal control problem.
\end{remark}

In the iterative process, $ H_{1}^{(i)} $ and $ H_{2}(i) $ are determined from the system data collected online, { which together with the vectors $ z_{k} $ form the quadratic type $Q_{1}^{(i)}, Q_{2}^{(i)}$}. In a deterministic system, we can directly separate the resulting $ H_{1} $ and $ H_{2} $ using the Kronecker product, i.e.,
$$ x^{T}Wy=(y^{T}\otimes x^{T})vec(W) $$
with vectors $ x\in \mathcal{R}^{vx} $ , $ y\in \mathcal{R}^{vy} $ and matrix $ W \in \mathcal{R}^{vx \times vy} $. After that, the least squares method is used to obtain the matrices $ H_{1} $ and $ H_{2} $ based on the collected data. In stochastic systems, a similar solution procedure is developed. The following equation is introduced to separate the matrices $ H_{1} $ and $ H_{2} $ from the quadratic Q-functions in order to facilitate the solution to their values.
$$ E(z^{T}_{k}Hz_{k}) = Tr\left[ E(z_{k}z_{k}^{T})H\right].  $$
Define {$ \mathcal{Z}_{k}$} as
$$ { \mathcal{Z}_{k}} := E(z_{k}z^{T}_{k})=\left[  \begin{matrix}	{z}_{11}&{z}_{12}&\cdots&{z}_{1p}\\{z}_{21}&{z}_{22}&\cdots&{z}_{1p}\\\vdots&\vdots&\ddots&\vdots\\{z}_{p1}&{z}_{p2}&\cdots&{z}_{pp}	  \end{matrix}\right]. $$
Then the two vectorization operators {are denoted as 	
		\begin{align}
			& vecs(H)=\left[
			\begin{matrix}H_{11}&H_{12}&H_{13}&\cdots&H_{1p}&H_{22}&\cdots&H_{pp}\end{matrix}\right]^{T},\nonumber\\
			& vech(\mathcal{Z}_{k})=\left[ \begin{matrix}
				{z}_{11}&2{z}_{12}&2{z}_{13}&\cdots&2{z}_{1{p}}&{z}_{22}&\cdots&{z}_{pp} \end{matrix}\right]\nonumber,
		\end{align}
where $ {vecs(H)} \in \mathcal{R}^{\frac{p(p+1)}{2}} $ and $ {vech(\mathcal{Z}_{k})}\in \mathcal{R}^{\frac{p(p+1)}{2}} $. $ H_{ij} $ represents the element in the $ i^{th} $ row and $ j^{th} $ column position of the symmetric matrix $ H $. Now $ Q^{*}(x_{k},u_{k},v_{k}) $ can be shown as a multiplication of two one-dimensional matrices
	$$ 	{Q^{*}(x_{k},u_{k},v_{k})=Tr\left[ E(z_{k}z_{k}^{T})H\right] =vech({\mathcal{Z}_{k}})vecs(H) }.$$
	Specifically, two Q-functions of the relevant performance during the iterative process are expressed as
	\begin{equation}\label{Qi+1:vector}
		\left\lbrace
		\begin{split}	
			&{Q_{1}^{(i+1)}(x_{k},u_{k}^{(i)},v_{k}^{(i)})=vech({\mathcal{Z}_{k}})vecs(H_{1}^{(i+1)})},\\
			&{Q_{2}^{(i+1)}(x_{k},u_{k}^{(i)},v_{k}^{(i)})=vech({\mathcal{Z}_{k}})vecs(H_{2}^{(i+1)})}.
		\end{split}\right.
	\end{equation}
	To solve $ Q^{(i+1)} $, the right-hand side of (\ref{Q1Q2:i}) is rewritten as
	\begin{equation}\label{dkH}
		\left\lbrace
		\begin{split}	
			&d_{1}(z_{k},H_{1}^{(i)})\\
			{:=}&r_{1}\left(  x_{k},u_{k}^{(i)},v_{k}^{(i)}\right)+Q_{1}^{(i)}(x_{k+1},u_{k+1}^{(i)},v_{k+1}^{(i)}),\\
			&d_{2}(z_{k},H_{2}^{(i)})\\
			{:=}&r_{2}\left(  x_{k},u_{k}^{(i)},v_{k}^{(i)}\right)+Q_{2}^{(i)}(x_{k+1},u_{k+1}^{(i)},v_{k+1}^{(i)}).
		\end{split}\right.
	\end{equation}
	The vector ${vech(\mathcal{Z}_{k})} $ is determined by the collected system state. At the $ (i+1)^{th} $ iteration, (\ref{dkH}) is well-established since $ H_{1}^{(i)} $ and $ H_{2}^{(i)} $ are known and the system state can be observed. The vector ${vecs({H}^{(i+1)})}  $, as an unknown variable, is obtained by minimizing the error in the least squares sense
	\begin{equation}\label{Ax=b}
		\left\lbrace
		\begin{split}	
			&	\min_{{vecs(H_{1}^{(i+1)})}}\left\| {vech({\mathcal{Z}_{k}})vecs(H_{1}^{(i+1)})}-d_{1}(z_{k},H_{1}^{(i)})\right\| ,\\
			&		\min_{{vecs(H_{2}^{(i+1)})}}\left\| {vech({\mathcal{Z}_{k}})vecs(H_{2}^{(i+1)})}-d_{2}(z_{k},H_{2}^{(i)})\right\| .
		\end{split}\right.
		\nonumber
	\end{equation}
	With a certain amount of system state information collected, we are able to determine $ {vecs({H}_{1}^{(i+1)})} $ and $ {vecs({H}_{2}^{(i+1)})} $ by (\ref{Ax=b}), and then find the control gain for the next update step.
	
	Note that both the control input $ u_{k} $ and external disturbance $ v_{k} $ in the vector $ z_{k} $ are in linear form with respect to the system state $ x_{k} $. Therefore, the matrix $ {vech(\mathcal{Z}_{k})} $ must not be of full rank. That means the equation (\ref{Ax=b}) in the least squares sense is not solvable. To overcome this problem, probing noises $ e_{uk}, e_{vk} $, that are also known as persistence of excitations, are added to $ u_{k} $ and $ v_{k} $, respectively. So $ u_{k}^{(i)} $ and $ v_{k}^{(i)} $ in the iteration are denoted as
	\begin{equation}
		\left\lbrace
		\begin{split}	
			&\hat{u}_{k}^{(i)}=K_{2}^{(i)}x_{k}+e_{uk},\\
			&\hat{v}_{k}^{(i)}=K_{1}^{(i)}x_{k}+e_{vk}.
		\end{split}\right.
	\end{equation}
	Then the key equations (\ref{Qi+1:vector}) and (\ref{dkH}) about the model-free algorithm are rewritten as
	\begin{equation}\label{eQi+1:vector}
		\left\lbrace
		\begin{split}	
			Q_{1}^{(i+1)}(x_{k},\hat{u}_{k},\hat{v}_{k})&=E(\hat{z}_{k}^{T}H^{(i+1)}_{1}\hat{z}_{k})\\
			&={vech({\hat{\mathcal{Z}}_{k}})vecs(H_{1}^{(i+1)})},\\
			Q_{2}^{(i+1)}(x_{k},\hat{u}_{k},\hat{v}_{k})&=E(\hat{z}_{k}^{T}H^{(i+1)}_{2}\hat{z}_{k})\\
			&={vech({\hat{\mathcal{Z}}_{k}})vecs(H_{2}^{(i+1)})},
		\end{split}\right.
	\end{equation}
	and
	\begin{equation}\label{edkH}
		\left\lbrace
		\begin{split}	
			&d_{1}(\hat{z}_{k},H_{1}^{(i)})=r_{1}\left(  x_{k},\hat{u}_{k}^{(i)},\hat{v}_{k}^{(i)}\right)+E(z_{k+1}^{T}H^{(i)}_{1}z_{k+1}|{\hat{z}_{k}}),\\
			&d_{2}(\hat{z}_{k},H_{2}^{(i)})=r_{2}\left(  x_{k},\hat{u}_{k}^{(i)},\hat{v}_{k}^{(i)}\right)+E(z_{k+1}^{T}H^{(i)}_{2}z_{k+1}|{\hat{z}_{k}}).
		\end{split}\right.
	\end{equation}
	
	{In order to calculate $H_{1}^{(i+1)}$ and $H_{2}^{(i+1)}$, we need to generate $N_{u}$ system trajectories on the time interval $[k,k+1]$ based on the state $x_{k}$, with the action of $u_{k}^{(i)}, v_{k}^{(i)}$. Then, the numerical averaging is used to approximate the expectation as
		\begin{equation}
			d(\hat{z}_{t},H^{(i)})=\dfrac{1}{N_{u}} \sum \limits_{j=0}^{N_{u}}\left[ r\left(  x_{k},\hat{u}_{k}^{(i)},\hat{v}_{k}^{(i)}\right)+z_{k+1}^{T}H^{(i)}z_{k+1}\right].
			\nonumber
	\end{equation}}
	After collecting $ N $ data tuples to ensure that the matrix $ X $ is of full rank, we obtain {$ {H}_{1}^{(i+1)},{H}_{2}^{(i+1)} $} as
	\begin{equation}\label{eAx=b}
		\left\lbrace
		\begin{split}	
			&{vecs(H_{1}^{(i+1)})}=(X^{T}X)^{-1}X^{T} Y_{1},\\
			&{vecs(H_{2}^{(i+1)})}=(X^{T}X)^{-1}X^{T} Y_{2},\\
		\end{split}\right.
	\end{equation}
	where
	\begin{equation}
		\begin{split}	
			&X=\left[ \begin{matrix}
				{vech({\hat{\mathcal{Z}}_{k}})^{T}}&{vech({\hat{\mathcal{Z}}_{k+1}})^{T}}&\cdots&{vech(\hat{\mathcal{Z}}_{k+N})^{T}}\end{matrix}\right] ^{T},\\
			&Y_{1}=\left[ \begin{matrix} d_{1}(\hat{z}_{k}, H_{1}^{(i)} \\ d_{1}(\hat{z}_{k+1},H_{1}^{(i)})\\\cdots \\ d_{1}(\hat{z}_{k+N},H_{1}^{(i)})\end{matrix} \right],Y_{2}=\left[ \begin{matrix}d_{2}(\hat{z}_{k}, H_{2}^{(i)}) \\ d_{2}(\hat{z}_{k+1},H_{2}^{(i)})\\ \cdots \\ d_{2}(\hat{z}_{k+N},H_{2}^{(i)})\end{matrix} \right] .
		\end{split} \nonumber
	\end{equation}
	In the calculation, we use the numerical average to approximate the expectations in equations (\ref{eQi+1:vector}) and (\ref{edkH}). The specific implementation steps of the model-free algorithm are presented as Algorithm 1.
	{In Algorithm 1, the convergence termination condition is set as
		\begin{itemize}
			\item[1)]  $ \left\| H_{1}^{(i+1)}-H_{1}^{(i)}\right\|<\varepsilon$.
			\item[2)]  $\left\| H_{2}^{(i+1)}-H_{2}^{(i)}\right\|<\varepsilon$.
			\item[3)]  $~~Q_{2}^{(i+1)}(x_{k},u^{(i+1)}_{k},v^{(i+1)}_{k})-Q_{2}^{(i)}(x_{k},u^{(i)}_{k},v_{k}^{(i)})$\\
			$<r_{2}(x_{k},u^{(i+1)}_{k},v^{(i+1)}_{k})$.
		\end{itemize}
		Conditions 1) and 2) are easily understood. In section 3.3, the reasons for the establishment of termination condition 3) are discussed in more details.}
	\begin{algorithm}
		\caption{{Procedure for Q-learning algorithm}.}
		\label{alg1}
		\begin{algorithmic}[1]
			\State \textbf{Input:} Number of data tuples $N$, tolerance $\varepsilon$, maximum number of iterations $i_{max}$, number of system trajectories $N_{u}$.
			
			\State \textbf{Initialization:} Initialize the iteration index $ i=0 $, the time index $ k=0 $. Set $ Q_{1}^{(0)}(\cdot)=0, Q_{2}^{(0)}(\cdot)=0 $. Choose any initial state vector $ x_{0} $ and the gains $ (K_{1}^{(0)},K_{2}^{(0)}) $.
			\For{$i=0$ to $i_{max }$}
			\State Applying $ \hat{u}_{t}^{(i)} $, $ \hat{v}_{t}^{(i)} $ with detection noises to system (\ref{system}).
			\State Calculate $ vech(\hat{\mathcal{Z}}_{k}), d_{1}(\hat{z}_{k},H_{1}^{(i)}), d_{2}(\hat{z}_{t},H_{2}^{(i)}). $
			\State Construct the matrices $X$, $Y_{1}$ and $Y_{2}$ on the time interval $[k, k+N] $.
			\State Estimate $vecs(H_{1}^{(i+1)})$, $vecs(H_{2}^{(i+1)})$ using \eqref{eAx=b}.
			\State Update $K_{1}^{(i+1)}$ and $K_{2}^{(i+1)}$ according to \eqref{K1K2:i by H1H2}.
			\If{{the iteration termination condition is satisfied} }
			\State break
			\EndIf
			\EndFor
			\State \textbf{Output:} The external disturbance and control strategy gains $K_{1}^{(i+1)},K_{2}^{(i+1)}$.
		\end{algorithmic}
	\end{algorithm}
	
	Considering that there are $ \frac{p(p+1)}{2} $ unknown elements in the vectors {$ vecs({H}_{1}^{(i+1)}), vecs({H}_{2}^{(i+1)})$} solved by least squares, the number of data tuples required should satisfy $ N \geq \frac{p(p+1)}{2} $.
	Given the initial values, $ H_{1}^{(i+1)} $ and $ H_{2}^{(i+1)} $ are obtained using (\ref{eAx=b}). Then they are brought into equation (\ref{K1K2:i by H1H2}) to calculate the gains $ (K_{1},K_{2}) $ for the next step. When $ i\rightarrow \infty $, $ Q^{(i)} \rightarrow Q^{*} $, $ (K^{(i)}_{1},K^{(i)}_{2})\rightarrow (K^{*}_{1},K^{*}_{2}) $. The convergence of this algorithm will be highlighted in Section 3.3.
	\begin{remark}
		The introduction of detection noises $ e_{uk},e_{vk} $ puts the system in a dynamic change so that we can gain more information about the system. Usually we set the persistence of excitation to be a white noise or trigonometric functions. However, it should be pointed that only the latter is permissible in the algorithm of this paper. Since the expectation of the former is a fixed value, {the utilization of detection noise in the form of a normal distribution is approximately equivalent to } adding a constant to some columns of the matrix $ X $. In this case, the reversibility of $ X $ does not change after collecting N data tuples.
	\end{remark}
\subsection{Convergence Analysis and Stability Proof}
	In this section, we will concentrate on analyzing the convergence of the the proposed Q-learning algorithm for the $ H_{2}/H_{\infty} $ problem and the stability of the final result. For ease of expression, we define
	\begin{equation}\label{AXY1Y2}
		\begin{split}
			&\mathcal{A}(X,Y_{1},Y_{2})\\
			:=&\left( A_{2}+C_{2}Y_{1}\right) ^{T}X\left( A_{2}+C_{2}Y_{1}\right) \\
			&+\left( A_{1}+B_{1}Y_{2}+C_{1}Y_{1}\right) ^{T}X\left( A_{1}+B_{1}Y_{2}+C_{1}Y_{1}\right) .
		\end{split}
	\end{equation}
	As a result, one has
	\begin{equation}\label{xk+1Pxk+1-xkPxk}
		\begin{split}
			E\left( x_{k+1}^{T}Px_{k+1}\right) =E\left(  x_{k}^{T}\mathcal{A}({P},{K_{1},K_{2}})x_{k}\right)
		\end{split}
	\end{equation}
	{with respect to system \eqref{system}.}
	\begin{lemma}\label{LemQ-P}
		The Q-learning algorithm based on equations (\ref{Q1Q2:i}) and (\ref{K1K2:i by H1H2}) is equivalent to the iterative relation of $ P_{1}^{(i)} $ and $ P_{2}^{(i)} $ as
		\begin{equation}
			\begin{split}
				&P_{1}^{(i+1)}=\left(A_{1}+B_{1}K_{2}^{(i+1)} \right)^{T}P_{1}^{(i)}\left(A_{1}+B_{1}K_{2}^{(i+1)} \right)-Q\\
				&+A_{2}^{T}P_{1}^{(i)}A_{2}-K_{2}^{(i+1)T}K_{2}^{(i+1)}-M_{1}^{(i+1)}\Delta_{1}^{-1}M_{1}^{(i+1)T},\\
				&P_{2}^{(i+1)}=\left(A_{1}+C_{1}K_{1}^{(i+1)} \right)^{T}P_{2}^{(i)}\left(A_{1}+C_{1}K_{1}^{(i+1)} \right)+Q\\
				&+\left(A_{2}+C_{2}K_{1}^{(i+1)}\right)^{T}P_{2}^{(i)}\left( A_{2}+C_{2}K_{1}^{(i+1)}\right)\\
				&-M_{2}^{(i+1)}\Delta_{2}^{-1}M_{2}^{(i+1)T},
			\end{split}
			\nonumber
		\end{equation}
		where
		\begin{equation}
			\begin{split}
				&M_{1}^{(i+1)}=\left(A_{1}+B_{1}K_{2}^{(i+1)} \right)^{T}P_{1}^{(i)}C_{1}+A_{2}^{T}P_{1}^{(i)}C_{2},\\
				&\Delta_{1}=\gamma^{2}I+C_{2}^{T}P_{1}^{(i)}C_{2}+C_{1}^{T}P_{1}^{(i)}C_{1},\\
				&M_{2}^{(i+1)}=\left(A_{1}+C_{1}K_{1}^{(i+1)} \right)^{T}P_{2}^{(i)}B_{1},\\
				&\Delta_{2}=I+B_{1}^{T}P_{2}^{(i)}B_{1}>0,
			\end{split}
			\nonumber
		\end{equation}
		and $$ K_{1}^{(i+1)}=-\Delta_{1}^{-1}M_{1}^{(i+1)T}, K_{2}^{(i+1)}=-\Delta_{2}^{-1}M_{2}^{(i+1)T} .$$
	\end{lemma}
{\bf Proof}.
		{See Appendix \ref{ProofLemQ-P}. $\blacksquare$}
	\begin{remark}
		In \cite{zhang2008infinite,zhang2007stochastic}, the authors proposed a numerical recursive algorithm similar to that in Lemma \ref{LemQ-P} for solving $ P_{1} $ and $ P_{2} $, and proved the convergence of the algorithm by the existence of the solution to Problem \ref{problem}. According to the results of \cite{zhang2008infinite,zhang2007stochastic}, combined with Lemma \ref{LemQ-P}, it is clear that the present algorithm in this paper is convergent. Unlike the proof in \cite{zhang2008infinite,zhang2007stochastic}, we prove the convergence of the sequence generated by the value function during the iterative processes, thus indicating that the mode-free RL algorithm in this article is convergent.
	\end{remark}
	
	Before proving the convergence and stability, inspired by the reference {\cite{rizvi2018output}}, we first give a theorem on persistence of excitation and prove that its presence does not lead to a bias in the results. Therefore, we can disregard the detection noise in the subsequent discussion, which simplifies the proof.
	\begin{theorem}\label{Thm-PEimpact}
		The introduction of the probing noise functions $ e_{uk} $ and $ e_{vk} $ does not cause any bias when estimating the matrices $ H_{1}^{(i)} $ and $ H_{2}^{(i)} $ in the Q-functions at each step.
	\end{theorem}
\textbf{{Proof}}.
		{According to Lemma \ref{LemQ-P}, the iterative process of Q-functions is essentially an iteration of $ (P_{1}^{(i)},P_{2}^{(i)}) $. Assume the estimated $ (H_{1}^{(i)},H_{2}^{(i)}) $ matrix after the introduction of noise is $ (\hat{H}_{1}^{(i)},\hat{H}_{2}^{(i)}) $ and the corresponding estimated $ (P_{1}^{(i)},P_{2}^{(i)}) $ is $ (\hat{P}_{1}^{(i)},\hat{P}_{2}^{(i)}) $. The Q-functions (\ref{eQi+1:vector}) can be rewritten as
		\begin{equation}
			\left\lbrace
			\begin{split}	
				Q_{1}^{(i+1)}&(x_{k},\hat{u}_{k},\hat{v}_{k})=E \begin{bmatrix}x_{k}\\\hat{u}_{k}^{(i)}\\\hat{v}_{k}^{(i)}\end{bmatrix}^{T}\hat{H}^{(i+1)}_{1}\begin{bmatrix}x_{k}\\\hat{u}_{k}^{(i)}\\\hat{v}_{k}^{(i)}\end{bmatrix} ,\\
				Q_{{ 2 }}^{(i+1)}&(x_{k},\hat{u}_{k},\hat{v}_{k})=E \begin{bmatrix}x_{k}\\\hat{u}_{k}^{(i)}\\\hat{v}_{k}^{(i)}\end{bmatrix}^{T}\hat{H}^{(i+1)}_{{ 2}}\begin{bmatrix}x_{k}\\\hat{u}_{k}^{(i)}\\\hat{v}_{k}^{(i)}\end{bmatrix},
			\end{split}\right.
			\nonumber
		\end{equation}
		and $ d_{1}(\hat{z}_{k},H_{1}^{(i)}), d_{2}(\hat{z}_{k},H_{2}^{(i)}) $ are changed as
		\begin{equation}
			\left\lbrace
			\begin{split}	
				&d_{1}(\hat{z}_{k},\hat{H}_{1}^{(i)})=r_{1}\left(  x_{k},\hat{u}_{k}^{(i)},\hat{v}_{k}^{(i)}\right)+E(z_{k+1}^{T}\hat{H}^{(i)}_{1}z_{k+1}),\\
				&d_{2}(\hat{z}_{k},\hat{H}_{2}^{(i)})=r_{2}\left(  x_{k},\hat{u}_{k}^{(i)},\hat{v}_{k}^{(i)}\right)+E(z_{k+1}^{T}\hat{H}^{(i)}_{2}z_{k+1}),
			\end{split}\right.
			\nonumber
		\end{equation}
where $ \hat{u}_{k}^{(i)}={u}_{k}^{(i)}+e_{uk} $ and $ \hat{v}_{k}^{(i)}={v}_{k}^{(i)}+e_{vk} $. We first demonstrate the unbiasedness of the estimate $ \hat{H}_{1}^{(i+1)} $. The Q-function $ Q_{1}^{(i+1)}(x_{k},\hat{u}_{k},\hat{v}_{k}) $ with noise $ \hat{u}_{k}^{(i)}, \hat{v}_{k}^{(i)} $ can be further expressed as
		\begin{equation}
			\begin{split}	
				&E \left( \begin{bmatrix}x_{k}\\{u}_{k}^{(i)}\\{v}_{k}^{(i)}\end{bmatrix}^{T}\hat{H}^{(i+1)}_{1}\begin{bmatrix}x_{k}\\ {u}_{k}^{(i)}\\{v}_{k}^{(i)}\end{bmatrix}+ \begin{bmatrix}0\\{e}_{uk}\\{e}_{vk}\end{bmatrix}^{T}\hat{H}^{(i+1)}_{1}\begin{bmatrix}x_{k}\\ {u}_{k}^{(i)}\\{v}_{k}^{(i)}\end{bmatrix}\right. \\
				&\left.~~~~+\begin{bmatrix}x_{k}\\\hat{u}_{k}^{(i)}\\{v}_{k}^{(i)}\end{bmatrix}^{T}\hat{H}^{(i+1)}_{1}\begin{bmatrix}0\\{e}_{uk}\\{e}_{vk}\end{bmatrix} +\begin{bmatrix}0\\{e}_{uk}\\{e}_{vk}\end{bmatrix}^{T}\hat{H}^{(i+1)}_{1}\begin{bmatrix}0\\{e}_{uk}\\{e}_{vk}\end{bmatrix} \right)\\
				&\overset{(\ref{Q1-Hi+1-Pi})}{=}E \left(  \begin{bmatrix}x_{k}^{T}&{u}_{k}^{(i)T}&{v}_{k}^{(i)T}\end{bmatrix}\hat{H}^{(i+1)}_{1}\begin{bmatrix}x_{k}^{T}&{u}_{k}^{(i)T}&{v}_{k}^{(i)T}\end{bmatrix}^{T}\right) \\
				&+E\left[  -e_{uk}^{T}u_{k}^{(i)}+\gamma^{2}e_{vk}^{T}v_{k}^{(i)}+\left( C_{2}e_{vk}\right)^{T}\hat{P}_{1}^{(i)}\left(A_{2}x_{k}+C_{2}v_{k}^{(i)} \right)\right. \\
				&+\left.  \left( B_{1}e_{uk}+C_{1}e_{vk}\right)^{T}\hat{P}_{1}^{(i)}\left( A_{1}x_{k}+B_{1}u_{k}^{(i)}+C_{1}v_{k}^{(i)}\right) \right]  \\
				&+E\left[ -u_{k}^{(i)T}e_{uk}+\gamma^{2}v_{k}^{(i)T}e_{vk}
				\right. \\
				&\left. +\left(A_{2}x_{k}+C_{2}v_{k}^{(i)} \right)^{T}\hat{P}_{1}^{(i)}\left(C_{2}e_{vk} \right)\right. \\
				&+\left.  \left( A_{1}x_{k}+B_{1}u_{k}^{(i)}+C_{1}v_{k}^{(i)}\right)^{T}\hat{P}_{1}^{(i)}\left( B_{1}e_{uk}+C_{1}e_{vk}\right) \right]  \\
				&+  E\left[ -e_{uk}^{T}e_{uk}+\gamma^{2}e_{vk}^{T}e_{vk}+\left(C_{2}e_{vk} \right)^{T}\hat{P}_{1}^{(i)}\left(C_{2}e_{vk} \right)\right. \\
				&+\left. \left( B_{1}e_{uk}+C_{1}e_{vk}\right)^{T}\hat{P}_{1}^{(i)}\left( B_{1}e_{uk}+C_{1}e_{vk}\right) \right].
			\end{split}
			\nonumber
		\end{equation}
		Next we deal with $ d_{1}(\hat{z}_{k},\hat{H}_{1}^{(i)}) $.
		According to system (\ref{system}), we easily obtain that
		\begin{equation}
			\begin{split}	
				&~~~~d_{1}(\hat{z}_{k},\hat{H}_{1}^{(i)})\\
				&=E \left(\gamma^{2}v_{k}^{(i)T}v_{k}^{(i)}-x_{k}^{T}Qx_{k}-u_{k}^{(i)T}u_{k}^{(i)}+\gamma^{2}v_{k}^{(i)T}e_{vk}\right. \\
				&\left.+\gamma^{2}e_{vk}^{T}v_{k}^{(i)}+\gamma^{2}e_{vk}^{T}e_{vk}-u_{k}^{(i)T}e_{uk} -e_{uk}^{T}u_{k}^{(i)}-e_{uk}^{T}e_{uk}\right) \\
				&+E\left[\left( A_{1}x_{k}+B_{1}{u}_{k}^{(i)} +C_{1}{v}_{k}^{(i)}\right)^{T}\right. \\
				&~~~~~~~~\left.\times\hat{P}_{1}^{(i)}\left( A_{1}x_{k}+B_{1}{u}_{k}^{(i)} +C_{1}{v}_{k}^{(i)}\right)\right. \\
				&\left. +\left( B_{1}{e}_{uk} +C_{1}{e}_{vk}\right)^{T}\hat{P}_{1}^{(i)}\left( A_{1}x_{k}+B_{1}{u}_{k}^{(i)} +C_{1}{v}_{k}^{(i)}\right)\right.    \\
				&\left. +\left( A_{1}x_{k}+B_{1}{u}_{k}^{(i)} +C_{1}{v}_{k}^{(i)}\right)^{T}\hat{P}_{1}^{(i)}\left( B_{1}{e}_{uk} +C_{1}{e}_{vk}\right) \right.    \\
				&\left. +\left(  B_{1}{e}_{uk} +C_{1}{e}_{vk}\right)^{T}\hat{P}_{1}^{(i)}\left( B_{1}{e}_{uk} +C_{1}{e}_{vk}\right) \right.   \\
				&\left. +  \left( A_{2}x_{k}+C_{2}{v}_{k}^{(i)}\right)^{T}\hat{P}_{1}^{(i)}\left( A_{2}x_{k}+C_{2}{v}_{k}^{(i)}\right)\right.\\
				&\left. +\left(C_{2}e_{vk} \right)^{T}\hat{P}_{1}^{(i)}\left(C_{2}e_{vk} \right)
				+  \left(C_{2}e_{vk} \right)^{T}\hat{P}_{1}^{(i)}\left(A_{2}x_{k}+C_{2}{v}_{k}^{(i)} \right) \right.\\
				& \left. +\left(A_{2}x_{k}+C_{2}{v}_{k}^{(i)} \right)^{T}\hat{P}_{1}^{(i)}\left(C_{2}e_{vk} \right)\right] .
			\end{split}
			\nonumber
		\end{equation}
		$ Q_{1}^{(i+1)}(x_{k},\hat{u}_{k},\hat{v}_{k}) $ and $ d_{1}(\hat{z}_{k},\hat{H}_{1}^{(i)}) $ are compared and the equation $ Q_{1}^{(i+1)}(x_{k},\hat{u}_{k},\hat{v}_{k})= d_{1}(\hat{z}_{k},\hat{H}_{1}^{(i)}) $ is handled by eliminating the same terms. After the treatment, there are no terms involving probing noise $ e_{uk} $ and $ e_{vk} $ on either side of the equation. Then
		\begin{equation}
			\begin{split}	
				&E \left(  \begin{bmatrix}x_{k}^{T}&{u}_{k}^{(i)T}&{v}_{k}^{(i)T}\end{bmatrix}\hat{H}^{(i+1)}_{1}\begin{bmatrix}x_{k}^{T}&{u}_{k}^{(i)T}&{v}_{k}^{(i)T}\end{bmatrix}^{T}\right)\\
				=&E \left[  \gamma^{2} {v}_{k}^{(i)T}{v}_{k}^{(i)}-x_{k}^{T}Qx_{k}-{u}_{k}^{(i)T}{u}_{k}^{(i)}\right] \\
				+&E\left[\left( A_{1}x_{k}+B_{1}{u}_{k}^{(i)} +C_{1}{v}_{k}^{(i)}\right)^{T}\hat{P}_{1}^{(i)}\right. \\
				&~~~~~~~~~~~~~~~~~~~~~~~~~~~~~~~~\left.\times \left( A_{1}x_{k}+B_{1}{u}_{k}^{(i)} +C_{1}{v}_{k}^{(i)}\right)\right. \\
				&\left. +  \left( A_{2}x_{k}+C_{2}{v}_{k}^{(i)}\right)^{T}\hat{P}_{1}^{(i)}\left( A_{2}x_{k}+C_{2}{v}_{k}^{(i)}\right)\right] .\\
			\end{split}
			\nonumber
		\end{equation}
		The above equation is identical to the first formula of (\ref{Q1Q2:i}) without the noises  $ e_{uk} $ and $ e_{vk} $, which shows $ \hat{H}_{1}^{(i+1)} = H_{1}^{(i+1)} $.  Based on (\ref{Q2-Hi+1-Pi}), the mechanism to prove $ \hat{H}_{2}^{(i+1)} = H_{2}^{(i+1)} $ is exactly the same as above, so proof process is delated.} $\blacksquare$
	\begin{lemma}\label{LemVI-P}
		The VI algorithm based on (\ref{Vi-VI-iteration}) and (\ref{viui-VI-iteration}), which is described briefly in Algorithm 2, is equivalent to the iteration relation of $ P_{1}^{(i)}, P_{2}^{(i)}$ as
		\begin{equation}
			\begin{split}
				&P_{1}^{(i+1)}=\left(A_{1}+B_{1}K_{2}^{(i)} \right)^{T}P_{1}^{(i)}\left(A_{1}+B_{1}K_{2}^{(i)} \right)-Q\\
				&+A_{2}^{T}P_{1}^{(i)}A_{2}-K_{2}^{(i)T}K_{2}^{(i)}-M_{1}^{(i)}\Delta_{1}^{-1}M_{1}^{(i)T},\\
				&P_{2}^{(i+1)}=\left(A_{1}+C_{1}K_{1}^{(i)} \right)^{T}P_{2}^{(i)}\left(A_{1}+C_{1}K_{1}^{(i)}x \right)+Q\\
				&+\left(A_{2}+C_{2}K_{1}^{(i)}\right)^{T}P_{2}^{(i)}\left( A_{2}+C_{2}K_{1}^{(i)}\right)-M_{2}^{(i)}\Delta_{2}^{-1}M_{2}^{(i)T},
			\end{split}
			\nonumber
		\end{equation}
		where
		\begin{equation}
			\begin{split}
				& M_{1}^{(i)}=\left(A_{1}+B_{1}K_{2}^{(i)} \right)^{T}P_{1}^{(i)}C_{1}+A_{2}^{T}P_{1}^{(i)}C_{2},\\
				&\Delta_{1}=\gamma^{2}I+C_{2}^{T}P_{1}^{(i)}C_{2}+C_{1}^{T}P_{1}^{(i)}C_{1}>0,\\
				& M_{2}^{(i)}=\left(A_{1}+C_{1}K_{1}^{(i)} \right)^{T}P_{2}^{(i)}B_{1},\\
				&\Delta_{2}=I+B_{1}^{T}P_{2}^{(i)}B_{1}>0,
			\end{split}
			\nonumber
		\end{equation}
		and $$ K_{1}^{(i)}=-\Delta_{1}^{-1}M_{1}^{(i)T},K_{2}^{(i)}=-\Delta_{2}^{-1}M_{2}^{(i)T} .$$
	\end{lemma}
	\begin{algorithm}
		\caption{{VI algorithm for $ H_{2}/H_{\infty} $ problem.}}
		\label{alg2}
		\begin{algorithmic}[1]
			\State \textbf{Initialization:} Start with $ P_{1}^{(0)}(\cdot)=0, P_{2}^{(0)}(\cdot)=0 $. Set the iteration index $ i=0 $ and the time index $ k=0 $.
			\State Obtain the initial control gains $ K_{1}^{(0)} = 0,  K_{2}^{(0)} = 0$ by (\ref{viui-VI-iteration}).
			\State \textbf{Value update}: 	
			\begin{equation}\label{Vi-VI-iteration}
				\begin{split}
					V_{1}^{(i+1)}(x_{k}) =&r_{1}\left(  x_{k},u_{k}^{(i)},v_{k}^{(i)}\right)+E\left( { x_{k+1}^{T} }P_{1}^{(i)}x_{k+1}\right),\\
					V_{2}^{(i+1)}(x_{k})=&r_{2}\left(  x_{k},u_{k}^{(i)},v_{k}^{(i)}\right)+E\left( { x_{k+1}^{T} }P_{2}^{(i)}x_{k+1}\right) .
				\end{split}
			\end{equation}
			\State \textbf{Policy improvement}:
			\begin{equation}\label{viui-VI-iteration}
				\begin{split}	
					v_{k}^{(i)}= \arg \min \limits_{v_{k}} &\left[ r_{1}\left(x_{k},u_{k}^{(i)},v_{k} \right) 	+E\left({ x_{k+1}^{T} }P_{1}^{(i)}x_{k+1}\right) \right] ,\\
					u_{k}^{(i)}	=\arg \min \limits_{u_{k}} &\left[  r_{2}\left(x_{k},u_{k},v_{k}^{(i)}\right) +E\left({ x_{k+1}^{T} }P_{2}^{(i)}x_{k+1}\right) \right].
				\end{split}
			\end{equation}
			\State $ i\leftarrow i+1 $ .
		\end{algorithmic}
	\end{algorithm}
\textbf{{Proof}}.
{Substituting $ M_{1}^{(i)},\Delta_{1},M_{2}^{(i)},\Delta_{2}, K_{1}^{(i)}$ and $K_{2}^{(i)} $ to the equation $ P_{1}^{(i+1)} $ and $ P_{2}^{(i+1)} $, it yields a form similar to (\ref{QlearnPiteration:A}).
\begin{equation}\label{VIPiteration:A}
			\left\lbrace \begin{split}
				P_{1}^{(i+1)}=&\mathcal{A}({P_{1}^{(i)}},{K_{1}^{(i)},K_{2}^{(i)}})-Q-K_{2}^{(i)T}K_{2}^{(i)}\\
				&+\gamma^{2}K_{1}^{(i)T}K_{1}^{(i)},\\
				P_{2}^{(i+1)}=&\mathcal{A}({P_{2}^{(i)}},{K_{1}^{(i)},K_{2}^{(i)}})+Q+K_{2}^{(i)T}K_{2}^{(i)}.
			\end{split}\right.
\end{equation}
{According to \eqref{xk+1Pxk+1-xkPxk}, the above equation \eqref{VIPiteration:A} is equivalent to (\ref{Vi-VI-iteration}) with any $ x_{k} $.} $\blacksquare$}
	
Equations (\ref{QlearnPiteration:A}) and (\ref{VIPiteration:A}) are exactly equivalent when $ (u_{k}, v_{k}) $ is substituted. Thus, the Q-learning algorithm and VI algorithm are also equivalent. Unlike the Q-function $ Q^{(i)}(x_{k},u_{k},v_{k}) $, $ V^{(i)}(x_{k}) $ involves only the system state $ x_{k} $, and its form is much simpler. Therefore, we show that the Q-learning algorithm is convergent by proving the convergence of the VI algorithm.
	
For a certain system state $ x_{k} $, the sequences generated during the iterative process are $ \left\lbrace V^{(i)}_{1}(x_{k})\right\rbrace $ and $ \left\lbrace V^{(i)}_{2}(x_{k}) \right\rbrace $. In general, we set the initial value $ Q^{(0)}(\cdot)=V^{(0)}(\cdot) = 0 $. For proving the convergence of the value iteration algorithm, some necessary lemmas are presented.
	\begin{lemma}\label{Lemfai>Vi}
		A new iterative formulation with initial values as $ \Psi_{1}^{(0)}(\cdot)=V^{(0)}_{1}(\cdot)=0, \Psi_{2}^{(0)}(\cdot)=V^{(0)}_{2}(\cdot)=0 $ is defined as
		\begin{equation}
			\left\lbrace \begin{split}
				&\Psi_{1}^{(i+1)}(x_{k}) = r_{1}\left(  x_{k},\eta_{2}^{(i)}x_{k},\eta_{1}^{(i)}x_{k}\right)+\Psi_{1}^{(i)}(x_{k+1}), \\
				&\Psi_{2}^{(i+1)}(x_{k}) = r_{2}\left(  x_{k},\eta_{2}^{(i)}x_{k},\eta_{1}^{(i)}x_{k}\right)+\Psi_{2}^{(i)}(x_{k+1}),
				\nonumber
			\end{split}\right.
		\end{equation}
		where the disturbance input $ \left\lbrace \eta_{1}^{(i)}x_{k}\right\rbrace$ and control input $ \left\lbrace \eta_{2}^{(i)}x_{k}\right\rbrace  $ are arbitrary. Then there exists a relationship between the sequences $ \left\lbrace \Psi^{(i)}(x_{k}) \right\rbrace $ and $ \left\lbrace V^{(i)}(x_{k})\right\rbrace$ as
		\begin{equation}\label{eq.fai>V}
			\begin{split}
				\Psi_{1}^{(i)}(x_{k}) \geq V_{1}^{(i)}(x_{k}), \Psi_{2}^{(i)}(x_{k}) \geq V_{2}^{(i)}(x_{k}).
			\end{split}
		\end{equation}
	\end{lemma}
\textbf{Proof}.
{According to (\ref{viui-VI-iteration}), it is obvious that  $ (u^{(i)},v^{(i)}) $ minimizes $ V_{1}^{(i+1)} $ and $ V_{2}^{(i+1)} $ at each iteration, while $ (\eta_{2}^{(i)}, \eta_{1}^{(i)})$ is arbitrary. So it is easy to obtain (\ref{eq.fai>V}). }$\blacksquare$
	
\begin{lemma}\label{LemV1V2dandiao}
		If the initial values are set to $ V_{1}^{(0)}(x_{k})=0 $ and $ V_{2}^{(0)}(x_{k})=0 $, then $ \left\lbrace V^{(i)}_{1}(x_{k})\right\rbrace  $ and  $ \left\lbrace V^{(i)}_{2}(x_{k})\right\rbrace  $ are nonincreasing and nondecreasing sequences, respectively, i.e., $ V_{1}^{(i+1)}(x_{k}) \leq V_{1}^{(i)}(x_{k})$,  $ V_{2}^{(i+1)}(x_{k})\geq V_{2}^{(i)}(x_{k}) $.
\end{lemma}
\textbf{Proof}.
		{See Appendix \ref{ProofLemV1V2dandiao}.} $\blacksquare$

	\begin{lemma} \label{LemV2youjie}
		There exists an upper bound $ Y(x_{k}) $ such that the sequence $ \left\lbrace V_{2}^{(i)}(x_{k}) \right\rbrace  $ satisfies $ 0\leq V_{2}^{(i)}(x_{k})\leq Y(x_{k}) $ .
	\end{lemma}
\textbf{Proof}.
		{See Appendix \ref{ProofLemV2youjie}.} $\blacksquare$
	
	\begin{lemma}\label{LemV1+V2>0}
		In the iterative process of VI algorithm, the value function at each step satisfies $ V_{1}^{(i)}(x_{k})+V_{2}^{(i)}(x_{k})\geq 0 $ with any state $ x_{k} $, i.e., $ P_{1}^{(i)}+P_{2}^{(i)}\geq 0 $ .
	\end{lemma}
\textbf{Proof}.
		{See Appendix \ref{ProofLemV1+V2>0}.} $\blacksquare$
	
	Next, we are ready to formally introduce the theorem on the convergence of Algorithm 1 as an essential conclusion based on Lemmas \ref{LemQ-P}-\ref{LemV1+V2>0}.
	\begin{theorem}\label{Thmconverge}
		As $ i \rightarrow \infty $, the sequences $ \left\lbrace u_{k}^{(i)} ,i \in \mathcal{N}_{+}\right\rbrace  $ and $ \left\lbrace v_{k}^{(i)} ,i \in \mathcal{N}_{+}\right\rbrace   $, that the Q-learning algorithm generates for solving Problem \ref{problem}, will converge to the optimal feedback control strategy $ u_{k}^{*} $ and the worst-case disturbance $v_{k}^{*}$, respectively, provided that the system (\ref{system}) is sufficiently excited.
	\end{theorem}
\textbf{Proof}.
{Both Q-learning and VI algorithm can essentially be viewed as the iterative relations for $ P^{(i)} $. Combining the conclusions of Lemma \ref{LemQ-P} and Lemma \ref{LemVI-P}, the Q-learning and VI algorithm are equivalent. In light of Lemma \ref{LemV2youjie}, it has been proven that $ 0\leq V_{2}^{(i+1)}(x_{k})\leq Y(x_{k}) $, i.e., $ -V_{2}^{(i+1)}(x_{k})\geq-Y(x_{k}) $. Combined with the conclusion of Lemma \ref{LemV1+V2>0}, one has $ V_{1}^{(i+1)}(x_{k})\geq -V_{2}^{(i+1)}(x_{k})\geq -Y(x_{k})$. Now we know that $ \left\lbrace V_{1}^{(i)}(x_{k})\right\rbrace  $ is a monotonically nonincreasing sequence with a lower bound; $\left\lbrace V_{2}^{(i)}(x_{k})\right\rbrace$ is a monotonically nondecreasing sequence with an upper bound. Therefore, during iterations, $ V_{1}^{(i)}(x_{k}) $ and $ V_{2}^{(i)}(x_{k}) $ converge to a unique value, respectively. Based on the equivalence of Algorithm 1 and Algorithm 2, it is clear that $ (Q_{1}^{(i)}, Q_{2}^{(i)}) \rightarrow (Q_{1}^{*},Q_{2}^{*}) $. As the iterative index $ i \rightarrow \infty $, there will have $ (K^{(i)}_{1},K^{(i)}_{2}) \rightarrow (K^{(\infty)}_{1}, K^{(\infty)}_{2}) = (K_{1}, K_{2}) $ given in (\ref{K1Riccati}) and (\ref{K2Riccati}).}  $\blacksquare$
	
	According to Theorem \ref{Thmconverge}, $ (u^{(i)}_{k},v^{(i)}_{k})\rightarrow (u_{k}^{*}, v_{k}^{*}) $ as $ i\rightarrow \infty $. The $ H_{2}/H_{\infty} $ controller $ u_{k}^{*} \in l_{\omega}^{2}(\mathcal{N}_{+}, \mathcal{R}^{m_{1}}) $ can make system (\ref{system}) ASMS with the {worst}-case disturbance $ v_{k}^{*} \in l_{\omega}^{2}(\mathcal{N}_{+}, \mathcal{R}^{m_{2}}) $. In other words, $ u_{k}^{(i)}$ is the admissible control after infinitely many iterations. Nevertheless, in the real-world application, an infinite number of iterations is impossible to achieve. Typically, the stopping criterion of VI algorithm is set as
	\begin{equation}
		\begin{split}
			&\left\| V_{1}^{(i+1)}(x_{k})-V_{1}^{(i)}(x_{k})\right\| < \varepsilon,\\
			& \left\| V_{2}^{(i+1)}(x_{k})-V_{2}^{(i)}(x_{k})\right\| < \varepsilon.\nonumber
		\end{split}
	\end{equation}
	One has $  V_{2}^{(i+1)}(x_{k+1})-V_{2}^{(i)}(x_{k+1}) < {\hat{\varepsilon}} $ based on Lemma \ref{Lemfai>Vi}. Consequently, when {the} break conditions are satisfied, the Lyapunov-type function is
	\begin{equation}
		V_{2}^{(i+1)}(x_{k+1})-V_{2}^{(i+1)}(x_{k}) < -r_{2}(x_{k},u^{(i)}_{k},v^{(i)}_{k})+{\hat{\varepsilon}}.
		\nonumber
	\end{equation}
	$ \varepsilon > 0 $ is an artificially specified tolerance error, which is determined by the actual demand. If the designed tolerance error satisfies $ {\hat{\varepsilon} }> r_{2}(x_{k},u^{(i)}_{k},v^{(i)}_{k}) $, the stability of the system will not be guaranteed under $ (u^{*}+\varepsilon_{1},v^{*}+\varepsilon_{2}) $. To overcome this difficulty, we add a more strict stopping condition to ensure that the obtained solution $ (u^{(i)}, v^{(i)}) $ to the $ H_{2}/H_{\infty} $ problem is capable of stabilizing the system.
	\begin{theorem}\label{Thmstable}
		For any $ x_{k}\neq 0 $, if the stopping condition
		\begin{equation}\label{stable-stopping}
			\begin{split}
				&Q_{2}^{(i+1)}(x_{k},u^{(i+1)}_{k},v^{(i+1)}_{k})-Q_{1}^{(i)}(x_{k},u^{(i)}_{k},v_{k}^{(i)})\\
				< & r_{2}(x_{k},u^{(i+1)}_{k},v^{(i+1)}_{k})
			\end{split}
		\end{equation}
		is satisfied, the pair of solutions $\left(u_{k}^{(i+1)}, v_{k}^{(i+1)}\right)$ for $ i \in \mathcal{N}_{k} $ stabilize system \eqref{system} asymptotically in the mean square sense.
	\end{theorem}
\textbf{Proof}.
{According to Lemma \ref{LemLyapunov}, in order to verify that  $(u^{(i+1)}, v^{(i+1)})$ can make the system ASMS, we need to determine whether
		\begin{equation}
			\begin{split}
				&\left( A_{1}+B_{1}K_{2}^{(i+1)}+C_{1}K_{1}^{(i+1)}\right) ^{T}P\left( A_{1}+B_{1}K_{2}^{(i+1)}\right. \\
				& \left. +C_{1}K_{1}^{(i+1)}\right) +(A_{2}+C_{2}K_{1}^{(i+1)})^{T}P(A_{2}+C_{2}K_{1}^{(i+1)})-P\\
				&=\mathcal{A}({P},{K_{1}^{(i+1)},K_{2}^{(i+1)}})-P < 0
			\end{split}
			\nonumber
		\end{equation}
		has a solution $ P>0 $. When equation (\ref{stable-stopping}) is met during the iteration, one has
		\begin{equation}
			\begin{split}
				&Q_{2}^{(i+1)}(x_{k},u^{(i+1)}_{k},v^{(i+1)}_{k})-Q_{1}^{(i)}(x_{k},u^{(i)}_{k},v_{k}^{(i)})\\
				&~~~~~~~~~~~~~~~~~~~~~~~~~~~~~~~~~~~~~-r_{2}(x_{k},u^{(i+1)}_{k},v^{(i+1)}_{k})\\
				=&E\left(x_{k}^{T} \begin{bmatrix}I\\K_{2}^{(i+1)}\\K_{1}^{(i+1)}	\end{bmatrix}^{T}H_{2}^{(i+1)}\begin{bmatrix}I\\K_{2}^{(i+1)}\\K_{1}^{(i+1)}	\end{bmatrix} x_{k}\right)
				\\
				&-E\left(x_{k}^{T} \begin{bmatrix}I\\K_{2}^{(i)}\\K_{1}^{(i)}	\end{bmatrix}^{T}H_{2}^{(i)}\begin{bmatrix}I\\K_{2}^{(i)}\\K_{1}^{(i)}	\end{bmatrix} x_{k}\right)\\
				&-E\left[ x_{k}^{T}(Q+K_{2}^{(i+1)T}K_{2}^{(i+1)})x_{k}\right] \\
				\overset{(\ref{Pi:Hi})}{=} &E(x_{k}^{T}P_{2}^{(i+1)}x_{k})-E[x_{k}^{T}\left( Q+K_{2}^{(i+1)T}K_{2}^{(i+1)}\right) x_{k}] \\
				&-E(x_{k}^{T}P_{2}^{(i)}x_{k})\\
				\overset{(\ref{QlearnPiteration:A})}{=}& E[x_{k}^{T}\mathcal{A}({P_{2}^{(i)}},{K_{1}^{(i+1)},K_{2}^{(i+1)}})x_{k}]-E(x_{k}^{T}P_{2}^{(i)}x_{k})\\
				<& 0.
			\end{split}
			\nonumber
		\end{equation}
		According to Lemma \ref{LemV1V2dandiao}, it is deduced that $ P^{(i)}_{2}>0 $ for $ i \in \mathcal{N}_{+} $. The Lyapunov-type inequality has at least one solution $ P_{2}^{(i)} > 0 $. As a consequence,  the solution $ (u^{(i)}, v^{(i)}) $ from (\ref{K1K2:i by H1H2}) can make system (\ref{system}) ASMS.}  $\blacksquare$
	\begin{remark}
		In general, the stopping criterion of the Q-learning algorithm is simply set to $\left\| H^{(i+1)}_{1}-H^{(i)}_{1}\right\| <\varepsilon$, $ \left\| H^{(i+1)}_{2}-H^{(i)}_{2}\right\| <\varepsilon $. Inspired by the research work \cite{wei2015value} on VI algorithm for deterministic systems, we add a stopping condition (\ref{stable-stopping}) for Algorithm 1. That ensures the control strategy $ u_{k}^{(i+1)} \in l_{\omega}^{2}(\mathcal{N}_{+}, \mathcal{R}^{m_{1}}) $ under the disturbance $ v_{k}^{(i+1)} \in l_{\omega}^{2}(\mathcal{N}_{+}, \mathcal{R}^{m_{2}}) $ is admissible.
	\end{remark}
	
	\section{ $ H_{2}/H_{\infty} $ controller design on F-16 aircraft autopilot}
	In this section, a simulation experiment on an F-16 aircraft autopilot model is provided to demonstrate that the Q-learning algorithm proposed in this paper can solve the stochastic $ H_{2}/H_{\infty} $ problem without using information about the dynamic model of the system.
	
	We improve an example of a deterministic discrete-time F-16 aircraft model given by {\cite{al2007model}}.  This model is assumed to suffer from multiplicative noise, since the system state $ x_{k} $ is always disturbed by the stochastic environment and this effect cannot be omitted. The F-16 aircraft model with random terms is described
	by \eqref{system} with matrices
	\begin{equation}
		\begin{split}
			&A_{1}=\left[\begin{matrix}0.906488&0.0816012&-0.0005\\0.0741349&0.90121& -0.000708383\\0&0&0.132655\end{matrix}\right]\\
			&B_{1}=\left[ \begin{matrix}-0.00150808\\-0.0096\\0.867345\end{matrix}\right], C_{1}=\left[ \begin{matrix}0.00951892\\ 0.00038373\\0\end{matrix}\right]\\
			&A_{2}=\left[ \begin{matrix}	0.0072&0.0026&0.0001\\0.0041&0.0917&0.0072\\0&0&0.0505\end{matrix}\right],C_{2}=\left[ \begin{matrix}0.00156\\0.00037\\0\end{matrix}\right].
		\end{split}
		\nonumber
	\end{equation}
	$ x_{k}=\left[ \begin{matrix}
		\alpha&q&\delta_{e}
	\end{matrix}  \right]^{T} $ represents the status information of the aircraft, including the angle of attack, pitch rate and elevator deflection angle. The disturbance attenuation is set to $ \gamma=1 $ and the parameter in the cost function can be known as $ Q=C^{T}C=I $. Employing the iterative algorithm in \cite{zhang2008infinite}, the coupling GAREs (\ref{P1Riccati}) and (\ref{P2Riccati}) associated with the stochastic F-16 aircraft system can be solved as
	\begin{equation}
		\begin{split}
			&P_{1}^{*}=\left[ \begin{matrix}-16.3448&-13.4481&0.0079\\-13.4481&-17.2342&0.0067\\0.0079&0.0067&-1.0101\end{matrix}\right],\\
			&P_{2}^{*}=\left[ \begin{matrix}16.9864&14.0870&-0.0082\\14.0870&17.8859&-0.0070\\-0.0082&-0.0070&1.0101\end{matrix}\right].
		\end{split}
		\nonumber
	\end{equation}
	The corresponding solution to Problem \ref{problem} is obtained as
	\begin{equation}
		\begin{split}
			&K_{1}^{*}=\left[ \begin{matrix}0.1559&0.1353&0\end{matrix}\right],\\
			&K_{2}^{*}=\left[ \begin{matrix}0.0949&0.1097&-0.0661\end{matrix}\right].
		\end{split}
		\nonumber
	\end{equation}
	
	Next, the Q-learning algorithm proposed in this paper is used to solve the stochastic mixed $ H_{2}/H_{\infty} $ problem for the third-order F-16 aircraft model. Suppose that the dynamics matrix is completely unknown. We initialize Algorithm 1 with $ H_{1}^{(0)}=0, H_{2}^{(0)}=0 $, the initial state $ x_{0}=\left[ \begin{matrix}10&5&-2\end{matrix} \right]^{T}$, and $ H_{2}/H_{\infty} $ gains $ K_{1}^{(0)}=\left[\begin{matrix}0.6305&1.6421&-1.0436\end{matrix}\right], K_{2}^{(0)}=\left[\begin{matrix}2.7695&0.1328&-0.1702\end{matrix}\right] $, that are not optimal.
	
	The iterative calculation is performed using the steps in Algorithm 1. For each iteration, 20 samples of data are collected to facilitate the estimation of matrices $ H_{1}^{(i)},H_{2}^{(i)} $ using the least squares method. The tolerance error is set to $ \varepsilon=10^{-3} $. {Referring to \cite{kiumarsi2017h}, the persistent excitation functions for the three cases  are given in order to place the system in dynamic change. After using the Q-learning algorithm proposed in this paper, the computational results $\hat{K}_{1}, \hat{K}_{2}$ are obtained for each case, respectively. And the variation of $P_{1}, P_{2}, K_{1}, K_{2}$, and system state are presented in Fig.\ref{Case1FigP}-Fig.\ref{Case3Figstate}.
		\begin{itemize}
			\item Case 1:
			\begin{equation}
				\begin{split}
					&e_{uk}=\sin(1.009k)+\cos^{2}(0.538k),\\
					&e_{vk}=\sin(9.7k)+\cos^{2}(10.2k).\\
					&\hat{K}_{1}=\left[ \begin{matrix}
						0.1556&0.1350&0
					\end{matrix}\right] ,\\
					&\hat{K}_{2}=\left[ \begin{matrix}
						0.0946&0.1093&-0.0661
					\end{matrix}\right] .
				\end{split}
				\nonumber
			\end{equation}
			\item Case 2:
			\begin{equation}
				\begin{split}
					&e_{uk}=\sin(0.9k)+\cos(100k),\\
					&e_{vk}=\sin(10k)+\cos(10k).\\
					&\hat{K}_{1}=\left[ \begin{matrix}
						0.1556&0.1350&0
					\end{matrix}\right] ,\\
					&\hat{K}_{2}=\left[ \begin{matrix}
						0.0946&0.1093&-0.0661
					\end{matrix}\right] .
				\end{split}
				\nonumber
			\end{equation}
			\item Case 3:
			\begin{equation}
				\begin{split}
					&e_{uk}= \sin(1.009k)+\cos^{2}(0.538k)+\sin(0.9k) \\
					&~~~~+\cos(100k),\\
					&e_{vk}=\sin(9.7k)+\cos^{2}(10.2k)+\sin(10k)+\cos(10k).\\
					&\hat{K}_{1}=\left[ \begin{matrix}
						0.1556&0.1350&0
					\end{matrix}\right] ,\\
					&\hat{K}_{2}=\left[ \begin{matrix}
						0.0946&0.1093&-0.0661
					\end{matrix}\right] .
				\end{split}
				\nonumber
			\end{equation}
	\end{itemize}}
	\begin{figure}[!htp]
		\centering
		\subfigure[$ \left\| P_{1}^{(i)}-P_{1}^{*}\right\| $]{\includegraphics[width=4.3cm,height=3.2cm]{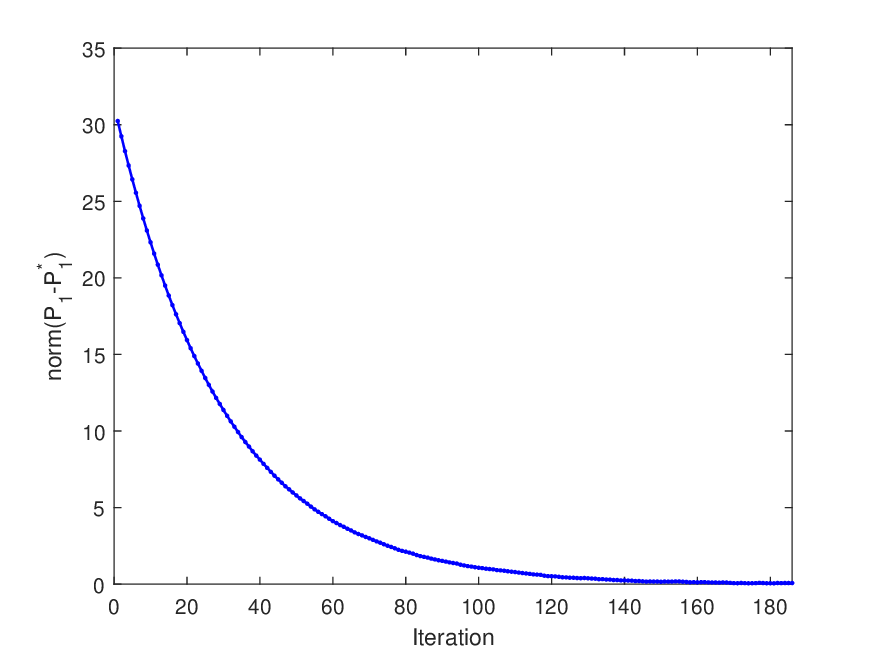}}
		\subfigure[$ \left\| P_{2}^{(i)}-P_{2}^{*}\right\| $]{\includegraphics[width=4.3cm,height=3.2cm]{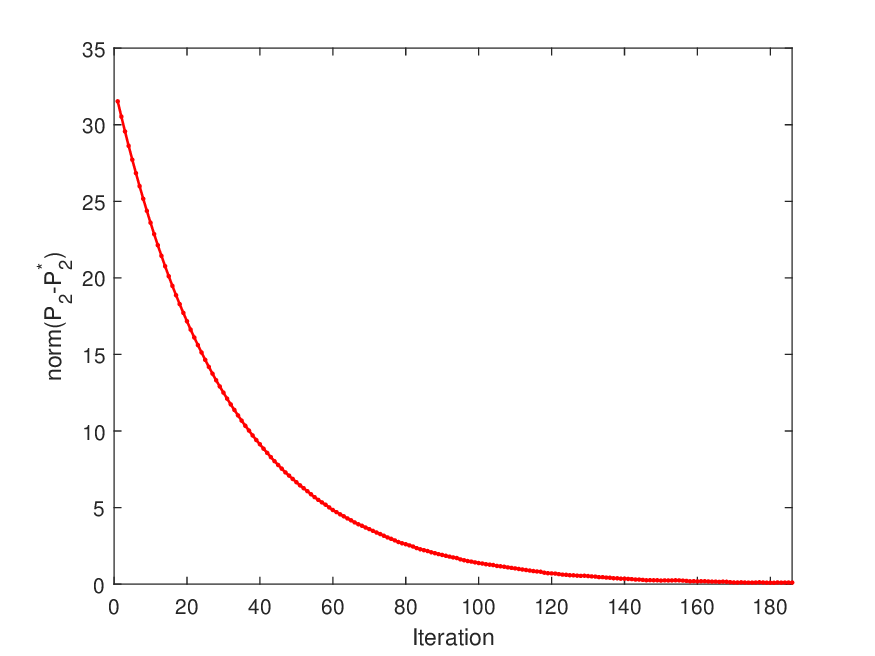}}
		\caption{Case 1: Convergence process of Q-function.}
		\label{Case1FigP}
	\end{figure}
	\begin{figure}[!htp]
		\centering
		\subfigure[$ \left\| P_{1}^{(i)}-P_{1}^{*}\right\| $]{\includegraphics[width=4.3cm,height=3.2cm]{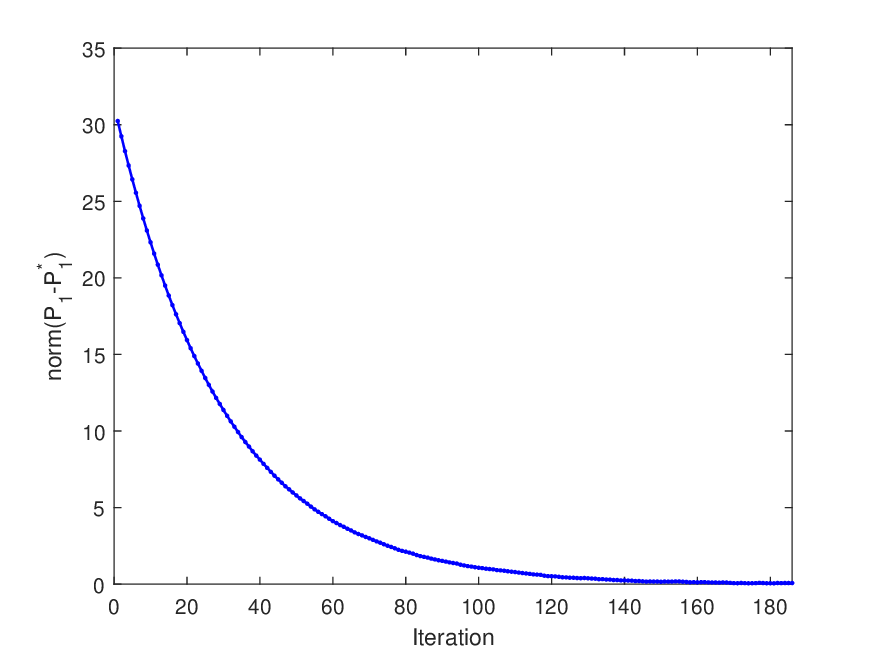}}
		\subfigure[$ \left\| P_{2}^{(i)}-P_{2}^{*}\right\| $]{\includegraphics[width=4.3cm,height=3.2cm]{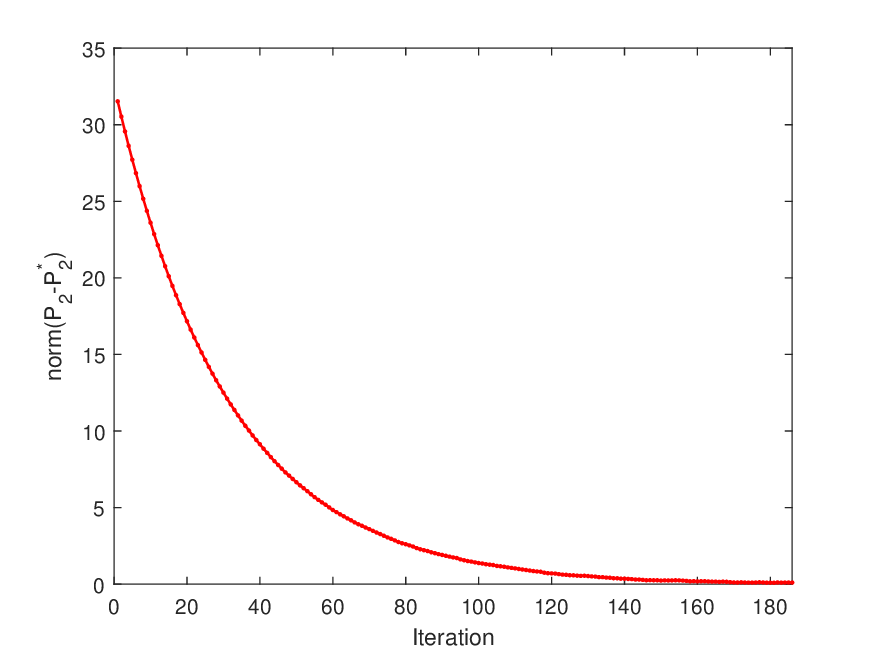}}
		\caption{Case 2: Convergence process of Q-function.}
		\label{Case2FigP}
	\end{figure}
	\begin{figure}[!htp]
		\centering
		\subfigure[$ \left\| P_{1}^{(i)}-P_{1}^{*}\right\|  $]{\includegraphics[width=4.3cm,height=3.2cm]{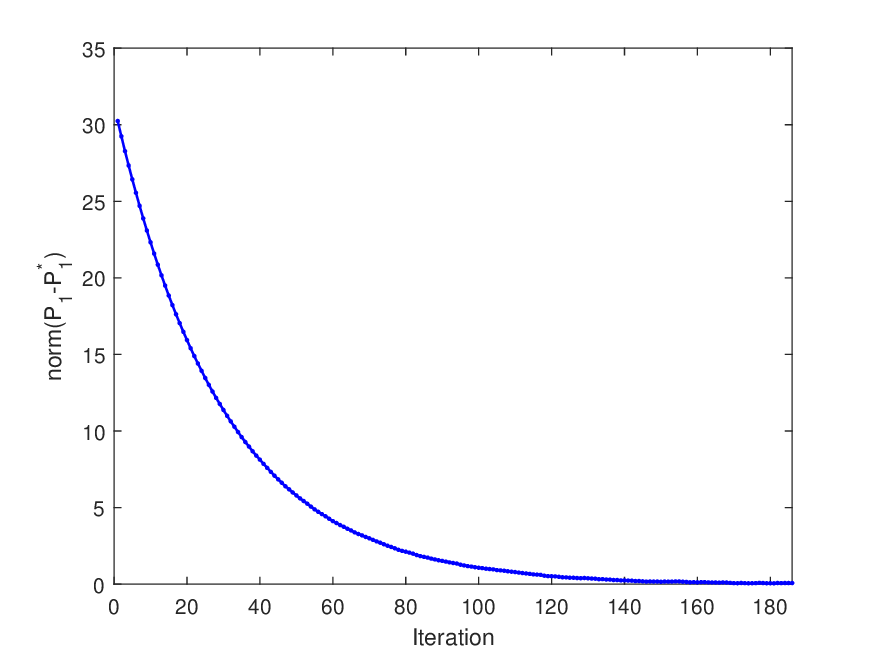}}
		\subfigure[$ \left\| P_{2}^{(i)}-P_{2}^{*}\right\|  $]{\includegraphics[width=4.3cm,height=3.2cm]{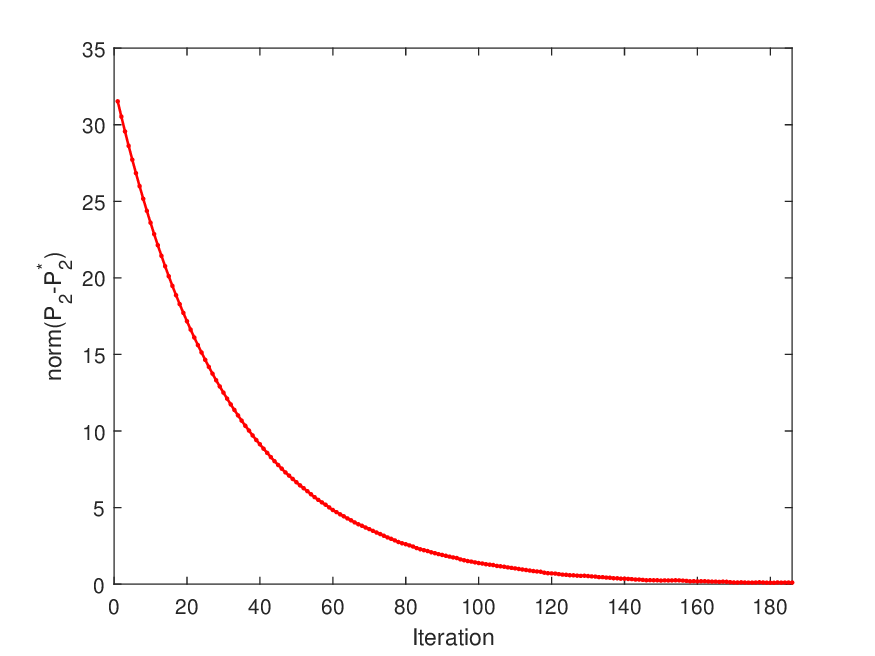}}
		\caption{Case 3: Convergence process of Q-function.}
		\label{Case3FigP}
	\end{figure}
	\begin{figure}[!htp]
		\centering
		\subfigure[$ \left\| K_{1}^{(i)}-K_{1}^{*}\right\|  $]{\includegraphics[width=4.3cm,height=3.2cm]{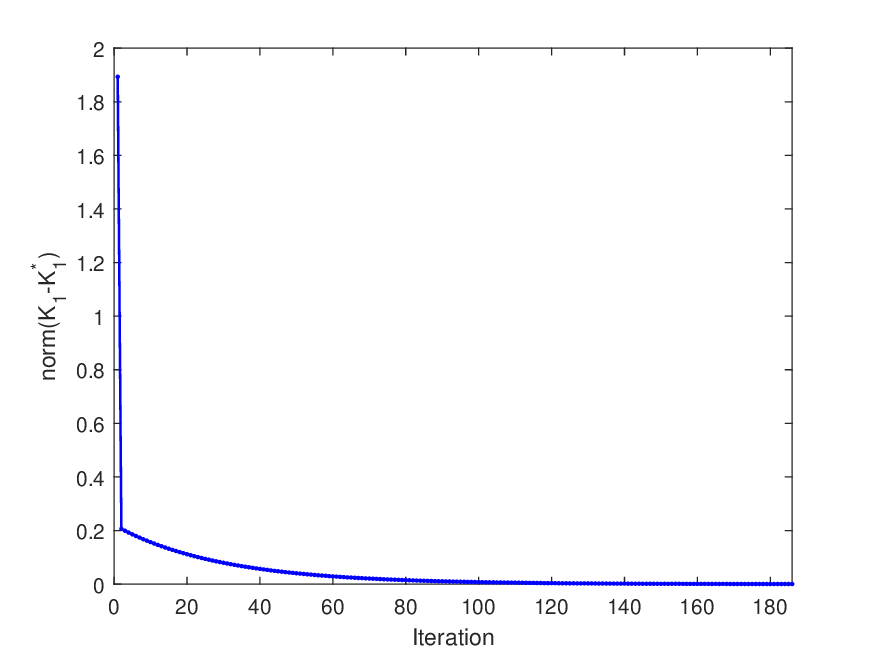}}
		\subfigure[$ \left\| K_{2}^{(i)}-K_{2}^{*}\right\|  $]{\includegraphics[width=4.3cm,height=3.2cm]{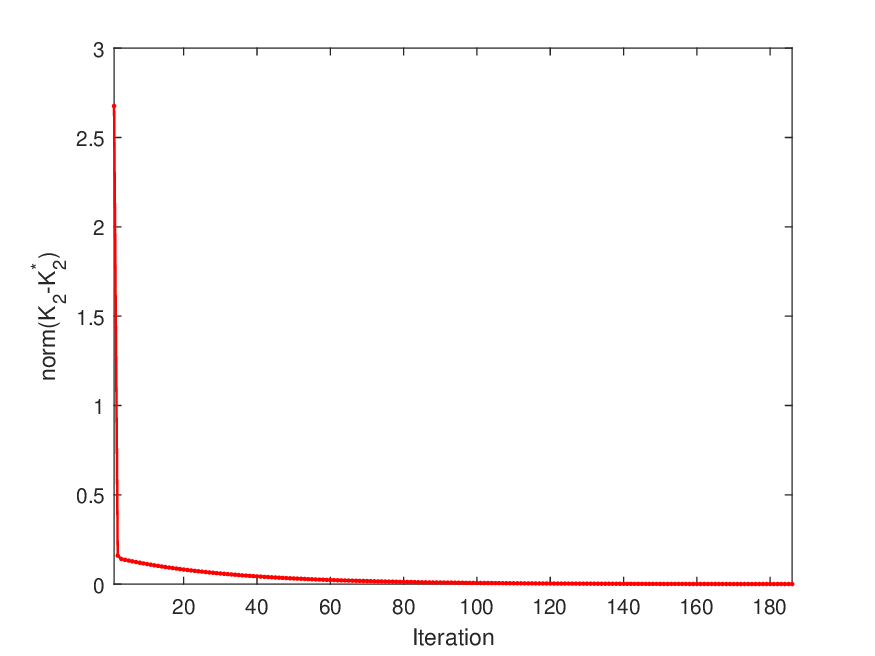}}
		\caption{Case 1: Convergence of $ K_{1} $ and $ K_{2} $ in Q-learning.}
		\label{Case1FigK}
	\end{figure}
	\begin{figure}[!htp]
		\centering
		\subfigure[$ \left\| K_{1}^{(i)}-K_{1}^{*}\right\|  $]{\includegraphics[width=4.3cm,height=3.2cm]{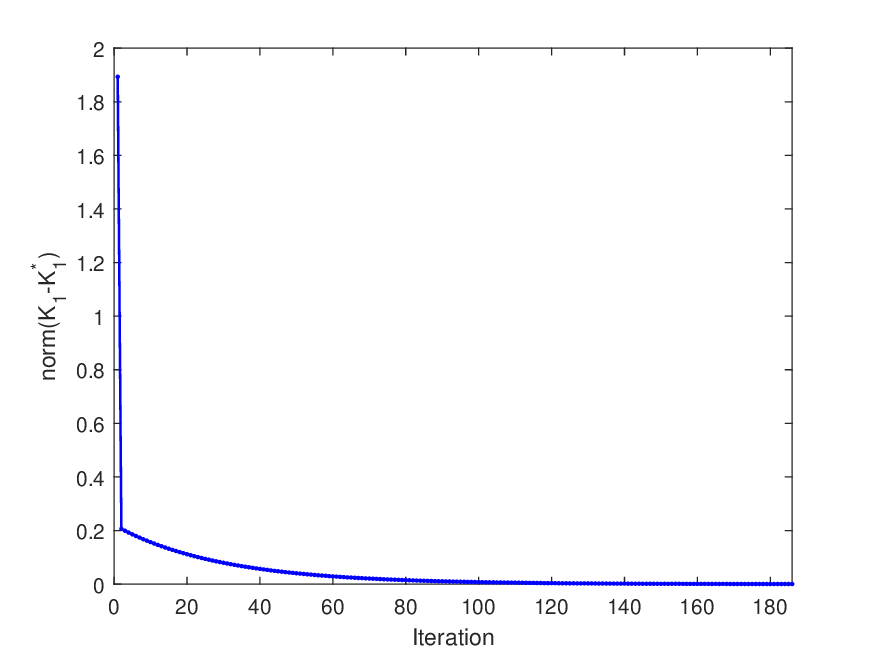}}
		\subfigure[$ \left\| K_{2}^{(i)}-K_{2}^{*}\right\|  $]{\includegraphics[width=4.3cm,height=3.2cm]{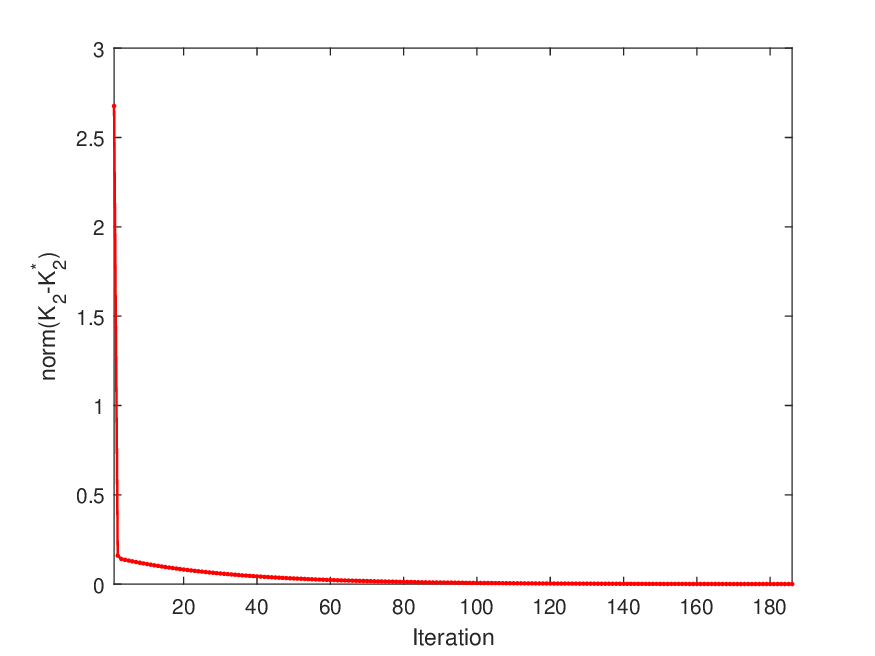}}
		\caption{Case 2: Convergence of $ K_{1} $ and $ K_{2} $ in Q-learning.}
		\label{Case2FigK}
	\end{figure}
	\begin{figure}[!htp]
		\centering
		\subfigure[$ \left\| K_{1}^{(i)}-K_{1}^{*}\right\|  $]{\includegraphics[width=4.3cm,height=3.2cm]{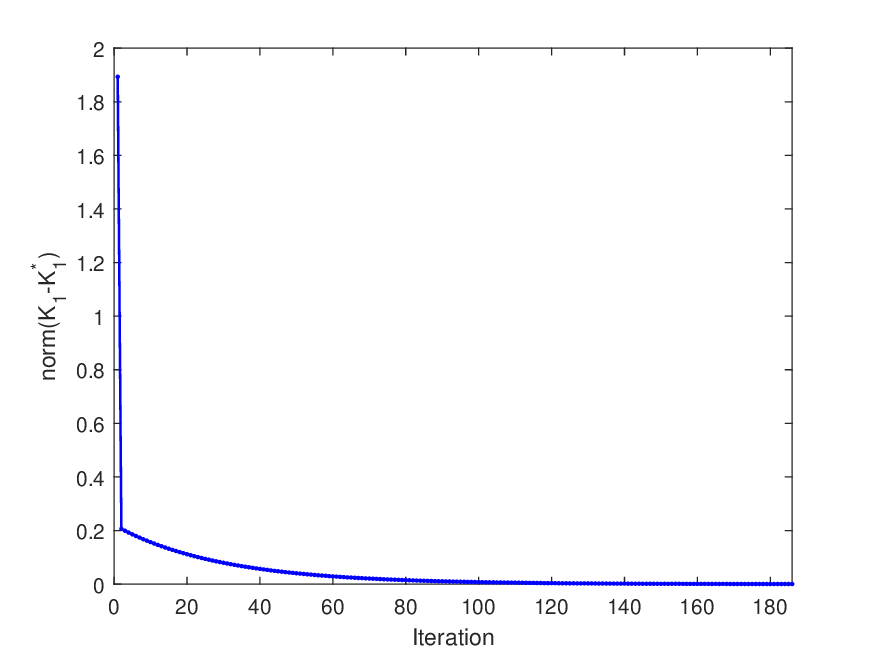}}
		\subfigure[$ \left\| K_{2}^{(i)}-K_{2}^{*}\right\|  $]{\includegraphics[width=4.3cm,height=3.2cm]{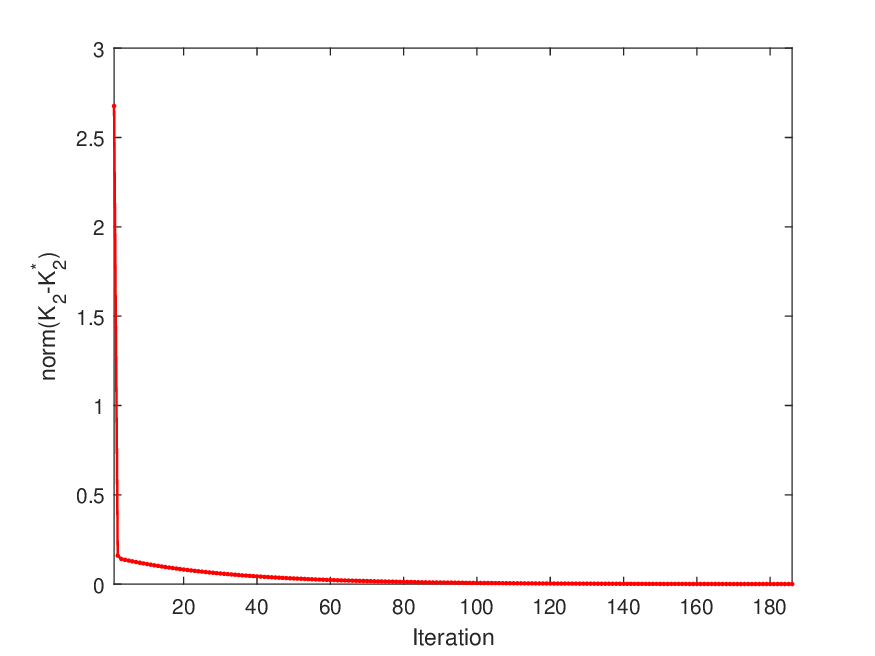}}
		\caption{Case 3: Convergence of $ K_{1} $ and $ K_{2} $ in Q-learning.}
		\label{Case3FigK}
	\end{figure}
	\begin{figure}[!htb]
		\centering
		\includegraphics[width=8cm,height=8cm]{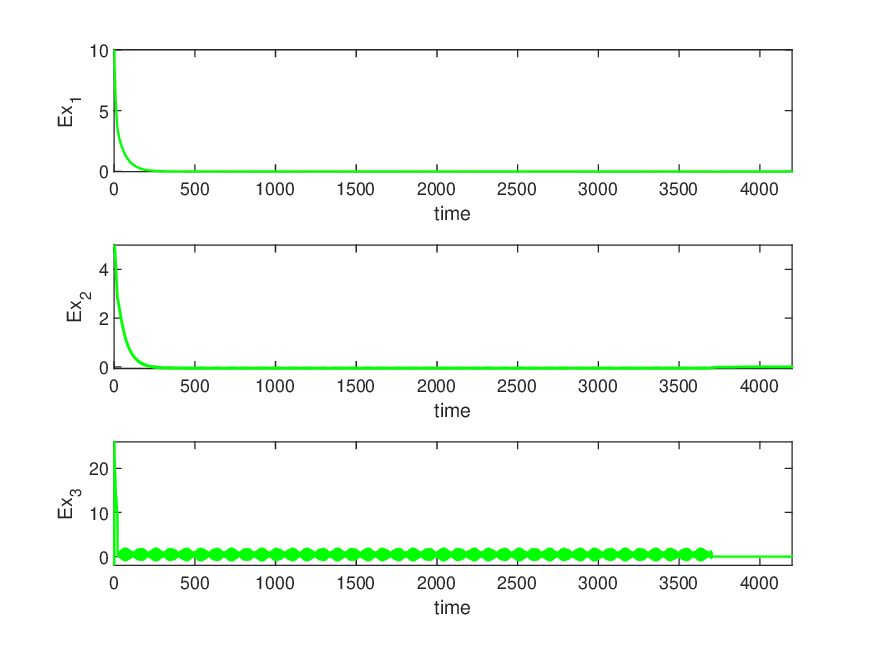}
		\caption{Case 1: The state curve in Q-learning algorithm.}
		\label{Case1Figstate}
	\end{figure}
	\begin{figure}[!htb]
		\centering
		\includegraphics[width=8cm,height=8cm]{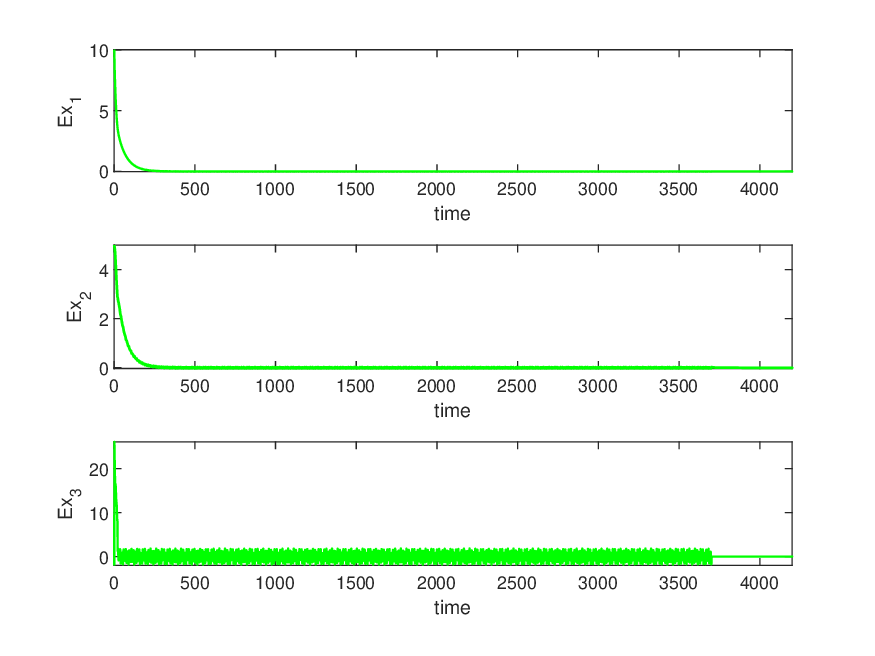}
		\caption{Case 2: The state curve in Q-learning algorithm.}
		\label{Case2Figstate}
	\end{figure}
	\begin{figure}[!htb]
		\centering
		\includegraphics[width=8cm,height=8cm]{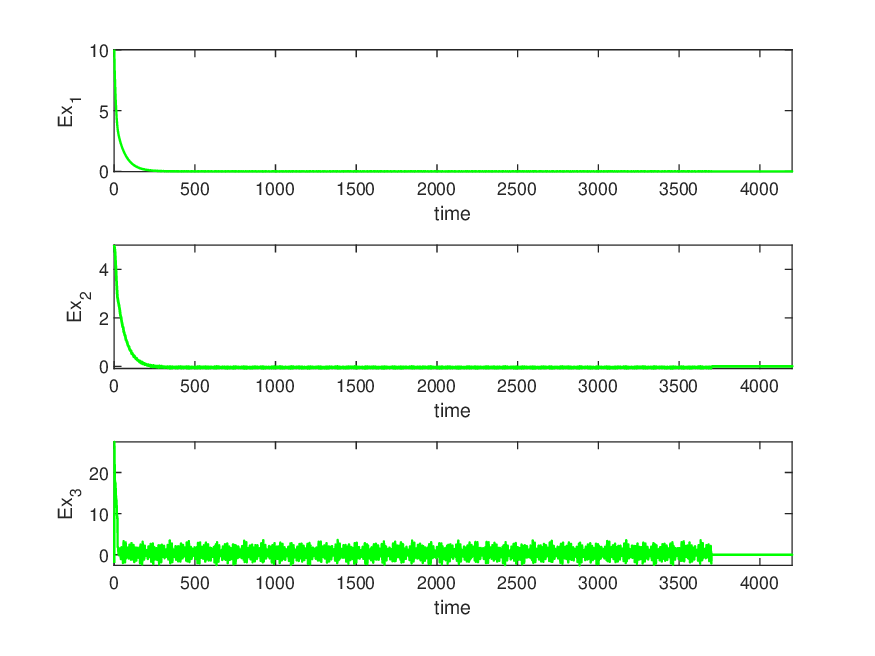}
		\caption{Case 3: The state curve in Q-learning algorithm.}
		\label{Case3Figstate}
	\end{figure}
{The Fig.\ref{Case1FigP}-Fig.\ref{Case3FigP} show how $ P_{1}^{(i)} $ and $ P_{2}^{(i)} $ change during the iterative process with three cases of probing noise} and demonstrates that the proposed algorithm is sufficient in convergence of $ (P_{1}^{(i)},P_{2}^{(i)}) $ to the optimal solution $ (P_{1}^{*},P_{2}^{*}) $. In {Fig.\ref{Case1FigK}-Fig.\ref{Case3FigK}}, we can know the variation of the corresponding gains $ K_{1}^{(i)} $ and $ K_{2}^{(i)} $. As can be seen from {Fig.\ref{Case1FigP}-Fig.\ref{Case3FigK}}, both the error between $ (P_{1}^{(i)},P_{2}^{(i)}) $ and $ (P_{1}^{*},P_{2}^{*}) $, and the error between $ (K_{1}^{(i)},K_{2}^{(i)}) $ and $ (K_{1}^{*},K_{2}^{*}) $, slowly become smaller until they converge to zero. Also we can notice from {Fig.\ref{Case1FigP}-Fig.\ref{Case3FigK}} that this deviation will remain fluctuating within a certain range when convergence is achieved. This is because of the bias that arises when the expectation is calculated using the average. When the amount of data increases in the process of calculating the average, the fluctuation range of this deviation becomes smaller.
	
	{Fig.\ref{Case1Figstate}-Fig.\ref{Case3Figstate} show} the change of the system state during the implementation of the algorithm, which is always affected by the detection noise. When the stopping criterion is satisfied, the detection noise terminates and the optimal control strategy obtained with the learning is sufficient to bring all the states to zero. Consequently, it can be shown that the obtained results $ u_{k} $ is the admissible $ H_{2}/H_{\infty} $ controller. Based on analyzing the results of the simulation, it can be concluded that the algorithm is effective.
	
	\section{Conclusion}
	In this paper, we have developed a model-free RL algorithm for stochastic systems with unknown system parameters to solve the designed problem. This algorithm can solve the optimal controller for the worst-case disturbance based on the system trajectory collected online. Compared with existing research works on Q-learning, the control policy can be admissible by means of a strict stopping criterion, and the unbiasedness introduced by the detection noise in the algorithm is demonstrated. Simulation experiments on an F-16 aircraft model with multiplicative noise have illustrated that the presented algorithm is effective. In our future research, it is a promising issue to reduce the number of data tuples used in the algorithm implementation as a way to improve efficiency.

\section{Appendix}
\subsection{Proof of Lemma \ref{LemQ-P}}\label{ProofLemQ-P}
Substituting  $ M_{1+1}^{(i)}, \Delta_{1}$ and $ K_{1}^{(i+1)} $ into the above iterative equation of $ P_{1}^{(i+1)} $ yields
\begin{equation}
	\begin{split}
		&P_{1}^{(i+1)}=\left(A_{1}+B_{1}K_{2}^{(i+1)} \right)^{T}P_{1}^{(i)}\left(A_{1}+B_{1}K_{2}^{(i+1)} \right)\\
		&+A_{2}^{T}P_{1}^{(i)}A_{2}+\left(A_{1}+B_{1}K_{2}^{(i+1)} \right)^{T}P_{1}^{(i)}C_{1}K_{1}^{(i+1)}\\
		&+A_{2}^{T}P_{1}^{(i)}C_{2}K_{1}^{(i+1)}-K_{2}^{(i+1)T}K_{2}^{(i+1)}-Q\\
		&=\left(A_{1}+B_{1}K_{2}^{(i+1)} +C_{1}K_{1}^{(i+1)}\right)^{T}P_{1}^{(i)}\left(A_{1}+B_{1}K_{2}^{(i+1)}\right. \\
		&+\left. C_{1}K_{1}^{(i+1)} \right)+A_{2}^{T}P_{1}^{(i)}\left( A_{2}+C_{2}K_{1}^{(i+1)}\right) \\
		&-\left(C_{1}K_{1}^{(i+1)}\right)^{T}P_{1}^{(i)}\left(A_{1}+B_{1}K_{2}^{(i+1)}+C_{1}K_{1}^{(i+1)} \right)\\
		&-K_{2}^{(i+1)T}K_{2}^{(i+1)}-Q\\
		&=\left(A_{1}+B_{1}K_{2}^{(i+1)} +C_{1}K_{1}^{(i+1)}\right)^{T}P_{1}^{(i)}\left(A_{1}+B_{1}K_{2}^{(i+1)}\right. \\
		&+\left. C_{1}K_{1}^{(i+1)} \right)+\left( A_{2}+C_{2}K_{1}^{(i+1)}\right) ^{T}P_{1}^{(i)}\left( A_{2}+C_{2}K_{1}^{(i+1)}\right) \\
		&-\left(C_{1}K_{1}^{(i+1)}\right)^{T}P_{1}^{(i)}\left(A_{1}+B_{1}K_{2}^{(i+1)}+C_{1}K_{1}^{(i+1)} \right)\\
		&-\left(C_{2}K_{1}^{(i+1)}\right)^{T}P_{1}^{(i)}\left(A_{2}+C_{2}K_{1}^{(i+1)}\right)\\
		&-K_{2}^{(i+1)T}K_{2}^{(i+1)}-Q\\
		&=\left(A_{1}+B_{1}K_{2}^{(i+1)} +C_{1}K_{1}^{(i+1)}\right)^{T}P_{1}^{(i)}\left(A_{1}+B_{1}K_{2}^{(i+1)}\right. \\
		&+\left.C_{1}K_{1}^{(i+1)} \right)+\left( A_{2}+C_{2}K_{1}^{(i+1)}\right) ^{T}P_{1}^{(i)}\left( A_{2}+C_{2}K_{1}^{(i+1)}\right)\\
		&-K_{1}^{(i+1)T}\left( M_{1}^{(i+1)T}+\Delta_{1}K_{1}^{(i+1)}-\gamma^{2}K_{1}^{(i+1)} \right)\\
		&-Q -K_{2}^{(i+1)T}K_{2}^{(i+1)} .
	\end{split}
	\nonumber
\end{equation}
Using a similar calculation,
there are
\begin{equation}
	\begin{split}
		P_{2}^{(i+1)}&=\left(A_{1}+C_{1}K_{1}^{(i+1)} \right)^{T}P_{2}^{(i)}\left(A_{1}+C_{1}K_{1}^{(i+1)} \right)\\
		&+\left(A_{2}+C_{2}K_{1}^{(i+1)}\right)^{T}P_{2}^{(i)}\left( A_{2}+C_{2}K_{1}^{(i+1)}\right)\\
		&+\left(A_{1}+C_{1}K_{1}^{(i+1)} \right)^{T}P_{2}^{(i)}B_{1}K_{2}^{(i+1)}+Q,
	\end{split}
	\nonumber
\end{equation}
and
\begin{equation}
	\begin{split}
		&~~~~P_{2}^{(i+1)}\\
		&=\left(A_{1}+B_{1}K_{2}^{(i+1)} +C_{1}K_{1}^{(i+1)}\right)^{T}P_{2}^{(i)}\left(A_{1}+B_{1}K_{2}^{(i+1)}\right. \\
		&+\left. C_{1}K_{1}^{(i+1)} \right)+\left(A_{2}+C_{2}K_{1}^{(i+1)}\right)^{T}P_{2}^{(i)}\left( A_{2}+C_{2}K_{1}^{(i+1)}\right)\\
		&-\left(B_{1}K_{2}^{(i+1)}\right)^{T}P_{2}^{(i)}\left(A_{1}+B_{1}K_{2}^{(i+1)}+C_{1}K_{1}^{(i+1)} \right)+Q\\
		&=\left(A_{1}+B_{1}K_{2}^{(i+1)} +C_{1}K_{1}^{(i+1)}\right)^{T}P_{2}^{(i)}\left(A_{1}+B_{1}K_{2}^{(i+1)}\right. \\
		&+\left. C_{1}K_{1}^{(i+1)} \right)-K_{2}^{(i+1)T}\left( M_{2}^{(i+1)T}+\Delta_{2}K_{2}^{(i+1)}-K_{2}^{(i+1)} \right) \\
		&+Q+\left( A_{2}+C_{2}K_{1}^{(i+1)}\right) ^{T}P_{1}^{(i)}\left( A_{2}+C_{2}K_{1}^{(i+1)}\right) .
	\end{split}
	\nonumber
\end{equation}
Based on (\ref{AXY1Y2}), one has
\begin{equation}\label{QlearnPiteration:A}
	\left\lbrace
	\begin{split}
		P_{1}^{(i+1)}=&\mathcal{A}({P_{1}^{(i)}},{K_{1}^{(i+1)},K_{2}^{(i+1)}})-Q\\
		&-K_{2}^{(i+1)T}K_{2}^{(i+1)}+\gamma^{2}K_{1}^{(i+1)T}K_{1}^{(i+1)},\\
		P_{2}^{(i+1)}=&\mathcal{A}({P_{2}^{(i)}},{K_{1}^{(i+1)},K_{2}^{(i+1)}})+Q\\
		&+K_{2}^{(i+1)T}K_{2}^{(i+1)}.
	\end{split} \right.
\end{equation}
We now handle the iterative formulation of previously proposed Q-learning algorithm. Changing the form of the first equation in (\ref{Q1Q2:i}), it can be rewritten as
\begin{equation}
	\begin{split}
		&~~~~E\begin{bmatrix}x_{k}^{T}&u_{k}^{(i)T}&v_{k}^{(i)T}	\end{bmatrix}H_{1}^{(i+1)} \begin{bmatrix}x_{k}^{T}&u_{k}^{(i)T}&v_{k}^{(i)T}	\end{bmatrix}^{T} \\
		&=E\left(\gamma^{2} v_{k}^{(i)T}v_{k}^{(i)}-x_{k}^{T}Qx_{k}-u_{k}^{(i)T}u_{k}^{(i)} \right)\\
		&+E\begin{bmatrix}x_{k+1}^{T}&u_{k+1}^{(i)T}&v_{k+1}^{(i)T}	\end{bmatrix}H_{1}^{(i)} \begin{bmatrix}x_{k+1}^{T}&u_{k+1}^{(i)T}&v_{k+1}^{(i)T}	\end{bmatrix}^{T} \\
		&= E\begin{bmatrix}x_{k}\\u_{k}^{(i)}\\v_{k}^{(i)}	\end{bmatrix}^{T}\begin{bmatrix}		-Q&0&0\\	0&-I&0\\	0&0&\gamma^{2}I	\end{bmatrix} \begin{bmatrix}x_{k}\\u_{k}^{(i)}\\v_{k}^{(i)}	\end{bmatrix}\\
		&+Ex_{k+1}^{T}\begin{bmatrix}I&K_{2}^{(i)T}&K_{1}^{(i)T}	\end{bmatrix}H_{1}^{(i)} \begin{bmatrix}I&K_{2}^{(i)T}&K_{1}^{(i)T}	\end{bmatrix}^{T}x_{k+1}\\
		&= E\begin{bmatrix}x_{k}\\u_{k}^{(i)}\\v_{k}^{(i)}	\end{bmatrix}^{T}\left\lbrace\begin{bmatrix}		-Q&0&0\\	0&-I&0\\	0&0&\gamma^{2}I	\end{bmatrix}\right. \\
		& \left. +\begin{bmatrix}A_{1}^{T}\\B_{1}^{T}\\C_{1}^{T}	\end{bmatrix}\begin{bmatrix}I\\K_{2}^{(i)}\\K_{1}^{(i)}	\end{bmatrix}^{T}H_{1}^{(i)} \begin{bmatrix}I\\K_{2}^{(i)}\\K_{1}^{(i)}	\end{bmatrix}\begin{bmatrix}A_{1}^{T}\\B_{1}^{T}\\C_{1}^{T}	\end{bmatrix}^{T}\right. \\
		&\left. +\begin{bmatrix}A_{2}^{T}\\0\\C_{2}^{T}	\end{bmatrix}\begin{bmatrix}I\\K_{2}^{(i)}\\K_{1}^{(i)}	\end{bmatrix}^{T}H_{1}^{(i)} \begin{bmatrix}I\\K_{2}^{(i)}\\K_{1}^{(i)}	\end{bmatrix}\begin{bmatrix}A_{2}^{T}\\0\\C_{2}^{T}	\end{bmatrix}^{T}\right\rbrace \begin{bmatrix}x_{k}\\u_{k}^{(i)}\\v_{k}^{(i)}	\end{bmatrix}.
	\end{split}
	\nonumber
\end{equation}
Considering {the system state  $ x_{k} $ is arbitrary}, it yields
\begin{equation}\label{Q1-Hi+1-Pi}
	\begin{split}
		H_{1}^{(i+1)}&=
		\begin{bmatrix}A_{2}^{T}\\0\\C_{2}^{T}	\end{bmatrix}P_{1}^{(i)}\begin{bmatrix}A_{2}^{T}\\0\\C_{2}^{T}	\end{bmatrix}^{T}\\
		&+ \begin{bmatrix}A_{1}^{T}\\B_{1}^{T}\\C_{1}^{T}	\end{bmatrix}P_{1}^{(i)}\begin{bmatrix}A_{1}^{T}\\B_{1}^{T}\\C_{1}^{T}	\end{bmatrix}^{T}+\begin{bmatrix}		-Q&0&0\\	0&-I&0\\	0&0&\gamma^{2}I	\end{bmatrix}.
	\end{split}
\end{equation}
Similarly, $ H_{2}^{(i+1)} $ can be given as
\begin{equation}\label{Q2-Hi+1-Pi}
	\begin{split}
		&H_{2}^{(i+1)}=\\
		&\begin{bmatrix}A_{1}^{T}\\B_{1}^{T}\\C_{1}^{T}	\end{bmatrix}P_{2}^{(i)}\begin{bmatrix}A_{1}^{T}\\B_{1}^{T}\\C_{1}^{T}	\end{bmatrix}^{T}+\begin{bmatrix}A_{2}^{T}\\0\\C_{2}^{T}	\end{bmatrix}P_{2}^{(i)}\begin{bmatrix}A_{2}^{T}\\0\\C_{2}^{T}	\end{bmatrix}^{T}+\begin{bmatrix}		Q&0&0\\	0&I&0\\	0&0&0	\end{bmatrix}.
	\end{split}
\end{equation}
According to equation (\ref{Pi:Hi}), we can further obtain
\begin{equation}\label{LammaP1i-Hi}
	\begin{split}
		&P_{1}^{(i+1)}= \begin{bmatrix}I&K_{2}^{(i+1)T}&K_{1}^{(i+1)T}	\end{bmatrix}H_{1}^{(i+1)}\begin{bmatrix}I\\K_{2}^{(i+1)}\\K_{1}^{(i+1)}	\end{bmatrix}\\
		=& \begin{bmatrix}I\\K_{2}^{(i+1)}\\K_{1}^{(i+1)}	\end{bmatrix}^{T} \begin{bmatrix}		-Q&0&0\\	0&-I&0\\	0&0&\gamma^{2}I	\end{bmatrix}\begin{bmatrix}I\\K_{2}^{(i+1)}\\K_{1}^{(i+1)}	\end{bmatrix}\\ &+\begin{bmatrix}I\\K_{2}^{(i+1)}\\K_{1}^{(i+1)}	\end{bmatrix}^{T}\begin{bmatrix}A_{1}^{T}\\B_{1}^{T}\\C_{1}^{T}	\end{bmatrix}P_{1}^{(i)}\begin{bmatrix}A_{1}^{T}\\B_{1}^{T}\\C_{1}^{T}	\end{bmatrix}^{T}\begin{bmatrix}I\\K_{2}^{(i+1)}\\K_{1}^{(i+1)}	\end{bmatrix}\\
		&+\begin{bmatrix}I\\K_{2}^{(i+1)}\\K_{1}^{(i+1)}\end{bmatrix}^{T}\begin{bmatrix}A_{2}^{T}\\0\\C_{2}^{T}	\end{bmatrix}P_{1}^{(i)}\begin{bmatrix}A_{2}^{T}\\0\\C_{2}^{T}	\end{bmatrix}^{T}\begin{bmatrix}I\\K_{2}^{(i+1)}\\K_{1}^{(i+1)}	\end{bmatrix}\\
		=&-{ Q }-K_{2}^{(i+1)T}K_{2}^{(i+1)}+\gamma^{2}K_{1}^{(i+1)T}K_{1}^{(i+1)}\\
		&+ (A_{1}+B_{1}K_{2}^{(i+1)} +C_{1}K_{1}^{(i+1)}) ^{T} \\
		&\times P_{1}^{(i)}(A_{1}+B_{1}K_{2}^{(i+1)}+C_{1}K_{1}^{(i+1)})  \\
		& +\left(  A_{2}+C_{2}K_{1}^{(i+1)}\right) ^{T}P_{1}^{(i+1)}\left(  A_{2}+C_{2}K_{1}^{(i+1)}\right) \\
		=&-Q-K_{2}^{(i+1)T}K_{2}^{(i+1)}+\gamma^{2}K_{1}^{(i+1)T}K_{1}^{(i+1)}\\
		&+\mathcal{A}({P_{1}^{(i)}},{K_{1}^{(i+1)},K_{2}^{(i+1)}}),
	\end{split}
\end{equation}
where $\mathcal{A}({P_{1}^{(i)}},{K_{1}^{(i+1)},K_{2}^{(i+1)}})$ and $\mathcal{A}({P_{2}^{(i)}},{K_{1}^{(i+1)},K_{2}^{(i+1)}}) $ are given as \eqref{AXY1Y2}.
In the same way, we can obtain $ P_{2}^{(i+1)} $ as
\begin{equation}\label{LammaP2i-Hi}
	\begin{split}
		&P_{2}^{(i+1)}
		= \begin{bmatrix}I&K_{2}^{(i+1)T}&K_{1}^{(i+1)T}	\end{bmatrix}H_{2}^{(i+1)}\begin{bmatrix}I\\K_{2}^{(i+1)}\\K_{1}^{(i+1)}	\end{bmatrix}\\
		&	=Q+K_{2}^{(i+1)T}K_{2}^{(i+1)}+\mathcal{A}({P_{2}^{(i)}},{K_{1}^{(i+1)},K_{2}^{(i+1)}}).
	\end{split}
\end{equation}
It is apparent that equations (\ref{LammaP1i-Hi})-(\ref{LammaP2i-Hi}) are equivalent to (\ref{QlearnPiteration:A}). $\blacksquare$
\subsection{Proof of Lemma \ref{LemV1V2dandiao}}\label{ProofLemV1V2dandiao}   
According to equations (\ref{Vi-VI-iteration}) and (\ref{viui-VI-iteration}), one has $$ P_{1}^{(0)}=0, P_{2}^{(0)}=0, K_{1}^{(0)}=0 ,K_{2}^{(0)}=0 $$ and $$V_{1}^{(1)}=-E\left( x_{k}^{T}Qx_{k}\right) , V_{2}^{(1)}=E\left( x_{k}^{T}Qx_{k}\right) .$$

(1) $ V_{1}^{(i+1)}(x_{k}) \leq V_{1}^{(i)}(x_{k}).$

This conclusion is proved by mathematical induction. When {$i=0$}, it is clear that $$ V_{1}^{(1)}(x_{k})= -E\left( x_{k}^{T}Qx_{k}\right) \leq V_{1}^{(0)}(x_{k})=0. $$ Assume that $ V_{1}^{(i+1)}(x_{k}) \leq V_{1}^{(i)}(x_{k})$ holds for $ i=n-1 $, so that $  P_{1}^{(n-1)}\geq P_{1}^{(n)} $ with any $ x_{k} $. Then for $ i=n $, we can know
\begin{equation}\label{V1:shuxueguina}
	\begin{split}
		&V_{1}^{(n+1)}(x_{k})\\
		=&r_{1}\left(x_{k},u_{k}^{(n)},v_{k}^{(n)}\right)+E\left[ x_{k+1}^{T}P_{1}^{(n)}x_{k+1}\right]\\
		=&r_{1}\left(x_{k},u_{k}^{(n)},v_{k}^{(n)}\right)+E\left[ x_{k}^{T}\mathcal{A}({P_{1}^{(n)}},{K_{1}^{(n)},K_{2}^{(n)}})x_{k}\right] .\\
	\end{split}
\end{equation}
According to Lemma \ref{Lemfai>Vi} at $\eta_{1}^{(i)}=K_{1}^{(n-1)}$ and $ \eta_{2}^{(i)}=K_{2}^{(n-1)} $ and $ P_{1}^{(n-1)}\geq P_{1}^{(n)}$, the equation (\ref{V1:shuxueguina}) can be further obtained that
\begin{equation}
	\begin{split}
		&V_{1}^{(n+1)}(x_{k})\\
		\leq& r_{1}\left(x_{k},u_{k}^{(n-1)},v_{k}^{(n-1)}\right)\\
		&+E[x_{k}^{T}\mathcal{A}({P_{1}^{(n)}},{K_{1}^{(n-1)},K_{2}^{(n-1)}})x_{k}]\\
		=& r_{1}\left(x_{k},u_{k}^{(n-1)},v_{k}^{(n-1)}\right)\\
		&+E\left[ x_{k}^{T}\mathcal{A}({P_{1}^{(n-1)}},{K_{1}^{(n-1)},K_{2}^{(n-1)}})x_{k}\right] \\
		&+E\left[ x_{k}^{T}\mathcal{A}({P_{1}^{(n)}-P_{1}^{(n-1)}},{K_{1}^{(n-1)},K_{2}^{(n-1)}})x_{k}\right] \\
		{\leq}& r_{1}\left(x_{k},u_{k}^{(n-1)},v_{k}^{(n-1)}\right)\\
		&+E[x_{k}^{T}\mathcal{A}({P_{1}^{(n-1)}},{K_{1}^{(n-1)},K_{2}^{(n-1)}})x_{k}]\\
		=&V_{1}^{(n)}(x_{k}).\\
	\end{split}
	\nonumber
\end{equation}
Hence, $ V_{1}^{(i+1)}(x_{k}) \leq V_{1}^{(i)}(x_{k})$ holds for any $ i \in \mathcal{N}_{k} $.

(2) $ V_{2}^{(i+1)}(x_{k}) \geq V_{2}^{(i)}(x_{k})$

The new sequence is defined as
$$ \widetilde{\Psi}_{2}^{(i+1)}(x_{k}) = r_{2}\left(  x_{k},\widetilde{\eta}_{2}^{(i)}x_{k},\widetilde{\eta}_{1}^{(i)}x_{k}\right)+\widetilde{\Psi}_{2}^{(i)}(x_{k+1}) $$
with $ \widetilde{\Psi}_{2}^{(0)}(\cdot)=0 $. According to Lemma \ref{Lemfai>Vi} at  $\widetilde{\eta}_{1}^{(i)}=K_{1}^{(i+1)}$ and $ \widetilde{\eta}_{2}^{(i)}=K_{2}^{(i+1)} $, it follows that $ \widetilde{\Psi}_{2}^{(i+1)}(x_{k}) \geq V_{2}^{(i+1)}(x_{k}) $.

In the next step, mathematical induction will be employed. For {$ i=0 $},
$$ V_{2}^{(1)}(x_{k})=E\left( x_{k}^{T}Qx_{k}\right) \geq \widetilde{\Psi}_{2}^{(0)}(x_{k})=0 .$$
Assume that
$$ V_{2}^{(i)}(x_{k})\geq \widetilde{\Psi}_{2}^{(i-1)}(x_{k}).$$
Since  $\widetilde{\eta}_{1}^{(i-1)}=K_{1}^{(i)}, \widetilde{\eta}_{2}^{(i-1)}=K_{2}^{(i)} $ are applied in system, it is clear that $ V_{2}^{(i)}(x_{k+1})\geq \widetilde{\Psi}_{2}^{(i-1)}(x_{k+1}) $ holds. According to the definitions of $ V_{2}^{(i+1)}(x_{k}) $ and $ \widetilde{\Psi}_{2}^{(i+1)}(x_{k}) $, there has
\begin{equation}
	\begin{split}
		&V_{2}^{(i+1)}(x_{k})-\widetilde{\Psi}_{2}^{(i)}(x_{k})\\
		=&r_{2}\left(x_{k},u_{k}^{(i)}\right.\left. ,v_{k}^{(i)}\right)+V_{2}^{(i)}(x_{k+1})\\
		&-r_{2}\left(x_{k},\widetilde{\eta}_{2}^{(i-1)}x_{k},\widetilde{\eta}_{1}^{(i-1)}x_{k}\right)-\widetilde{\Psi}_{2}^{(i-1)}(x_{k+1})\\
		=&V_{2}^{(i)}(x_{k+1})-\widetilde{\Psi}_{2}^{(i-1)}(x_{k+1})	\geq 0.
	\end{split}
	\nonumber
\end{equation}
Hence it is determined that $ V_{2}^{(i+1)}(x_{k})\geq\widetilde{\Psi}_{2}^{(i)}(x_{k}) $ is constant. Combined with Lemma \ref{Lemfai>Vi}, it follows that
$$  V_{2}^{(i+1)}(x_{k})\geq \widetilde{\Psi}_{2}^{(i)}(x_{k})\geq V_{2}^{(i)}(x_{k}),$$
i.e., $  V_{2}^{(i+1)}(x_{k})\geq V_{2}^{(i)}(x_{k})$. $\blacksquare$
\subsection{Proof of Lemma \ref{LemV2youjie}} \label{ProofLemV2youjie}
Set $ u_{k}=\eta_{2}x_{k} $ as any determined admissible control under the disturbance $ v_{k}=\eta_{1}x_{k} $, and define the sequence update to satisfy
$$ \hat{\Psi}_{2}^{(i+1)}(x_{k}) = r_{2}\left(  x_{k},\eta_{2}x_{k},\eta_{1}x_{k}\right)+\hat{\Psi}_{2}^{(i)}(x_{k+1}),$$
where the initial value is $ \hat{\Psi}_{2}^{(0)}(\cdot)=0 $. So
\begin{equation}\label{fai-faiditui}
	\begin{split}
		&\hat{\Psi}_{2}^{(i+1)}(x_{k})-\hat{\Psi}_{2}^{(i)}(x_{k})\\
		=&r_{2}\left(x_{k},\eta_{2}x_{k},\eta_{1}x_{k}\right)+\hat{\Psi}_{2}^{(i)}(x_{k+1})\\
		&-r_{2}\left(x_{k},\eta_{2}x_{k},\eta_{1}x_{k}\right)-\hat{\Psi}_{2}^{(i-1)}(x_{k+1})\\
		=&\hat{\Psi}_{2}^{(i)}(x_{k+1})-\hat{\Psi}_{2}^{(i-1)}(x_{k+1}).
	\end{split}
\end{equation}
The update strategies for each step are $ \eta_{1}x_{k}, \eta_{2}x_{k} $. According to (\ref{fai-faiditui}), it can be deduced that
\begin{equation}
	\begin{split}
		\hat{\Psi}_{2}^{(i+1)}(x_{k})-\hat{\Psi}_{2}^{(i)}(x_{k})&=\hat{\Psi}_{2}^{(i)}(x_{k+1})-\hat{\Psi}_{2}^{(i-1)}(x_{k+1})\\
		&=\hat{\Psi}_{2}^{(i-1)}(x_{k+2})-\hat{\Psi}_{2}^{(i-2)}(x_{k+2})\\
		&~~~~\vdots\\
		&=\hat{\Psi}_{2}^{(1)}(x_{k+i})-\hat{\Psi}_{2}^{(0)}(x_{k+i})\\
		&	=\hat{\Psi}_{2}^{(1)}(x_{k+i}).
	\end{split}
	\nonumber
\end{equation}
From the above equation, it can be further derived that
\begin{equation}
	\begin{split}
		&\hat{\Psi}_{2}^{(i+1)}(x_{k})=\hat{\Psi}_{2}^{(1)}(x_{k+i})+\hat{\Psi}_{2}^{(i)}(x_{k})\\
		=&\hat{\Psi}_{2}^{(1)}(x_{k+i})+\hat{\Psi}_{2}^{(1)}(x_{k+i-1})+\hat{\Psi}_{2}^{(i-1)}(x_{k})\\
		=&\hat{\Psi}_{2}^{(1)}(x_{k+i})+\hat{\Psi}_{2}^{(1)}(x_{k+i-1})+\hat{\Psi}_{2}^{(1)}(x_{k+i-2})\\
		&+\hat{\Psi}_{2}^{(i-2)}(x_{k})\\
		&\vdots\\
		=&\hat{\Psi}_{2}^{(1)}(x_{k+i})+\hat{\Psi}_{2}^{(1)}(x_{k+i-1})+\hat{\Psi}_{2}^{(1)}(x_{k+i-2})\\
		&+{\hat{\Psi}_{2}^{(1)}(x_{k+i-3})}+\cdots+\hat{\Psi}_{2}^{(1)}(x_{k})+\hat{\Psi}_{2}^{(0)}(x_{k})\\
		=&\sum_{n=0}^{i}\hat{\Psi}_{2}^{(1)}(x_{k+n})\\
		=&\sum_{n=0}^{i} E\left[ x_{k+n}^{T}\left( Q+\eta_{2}^{T}\eta_{2}\right) x_{k+n}\right]
	\end{split}\nonumber\
\end{equation}
with $ \hat{\Psi}_{2}^{(0)}(x_{k})=0 $. Based on the conclusion of Lemma \ref{Lemfai>Vi}, it is known that
\begin{equation}
	\begin{split}
		V_{2}^{(i+1)}(x_{k})&\leq\hat{\Psi}_{2}^{(i+1)}(x_{k})\\
		&=\sum_{n=0}^{i} E\left[ x_{k+n}^{T}\left( Q+\eta_{2}^{T}\eta_{2}\right) x_{k+n}\right] \\
		&\leq \sum_{n=0}^{\infty} E\left[ x_{k+n}^{T}\left( Q+\eta_{2}^{T}\eta_{2}\right) x_{k+n}\right] .
	\end{split}
	\nonumber
\end{equation}
Due to the admissibility of the control input, there exists an upper bound $ Y(x_{k}) $ to satisfy that
\begin{equation}
	\begin{split}
		0\leq V_{2}^{(i+1)}(x_{k}) & \leq \sum_{n=0}^{\infty} E\left[ x_{k+n}^{T}\left( Q+\eta_{2}^{T}\eta_{2}\right) x_{k+n}\right]\\
		& \leq Y(x_{k}).
	\end{split}
	\nonumber
\end{equation} $\blacksquare$
\subsection{Proof of Lemma \ref{LemV1+V2>0}}\label{ProofLemV1+V2>0}
According to the iterative equations (\ref{Vi-VI-iteration}) and (\ref{viui-VI-iteration}) of the VI algorithm, it is known that $ P_{1}^{(0)}=0, P_{2}^{(0)}=0, K_{1}^{(0)}=0, K_{2}^{(0)}=0 $ and $ P_{1}^{(1)}=-Q, P_{2}^{(1)}=Q $. Obviously, $$ V_{1}^{(0)}(x_{k})+V_{2}^{(0)}(x_{k})=0, V_{1}^{(1)}(x_{k})+V_{2}^{(1)}(x_{k})=0$$ satisfies $ V_{1}^{(i)}(x_{k})+V_{2}^{(i)}(x_{k})\geq0 $ for $ i=0,1 $.

Suppose $ V_{1}^{(i)}(x_{k})+V_{2}^{(i)}(x_{k})\geq 0 $, according to the iterative formula, it can be seen that
\begin{equation}
	\begin{split}
		&V_{1}^{(i+1)}(x_{k})+V_{2}^{(i+1)}(x_{k})\\
		=&r_{1}\left(  x_{k},u_{k}^{(i)},v_{k}^{(i)}\right)+{V_{1}^{(i)}(x_{k+1})}\\
		&+r_{2}\left(  x_{k},u_{k}^{(i)},v_{k}^{(i)}\right)+{V_{2}^{(i)}(x_{k+1})}\\
		=&E\left[ \gamma^{2}x_{k}^{T}K_{1}^{(i)T}K_{1}^{(i)}x_{k}\right] +{V_{1}^{(i)}(x_{k+1})+V_{2}^{(i)}(x_{k+1})}\\
		\geq& V_{1}^{(i)}(x_{k+1})+V_{2}^{(i)}(x_{k+1})\\
		\geq& 0.
	\end{split}
	\nonumber
\end{equation}
Therefore, $ V_{1}^{(i)}(x_{k})+V_{2}^{(i)}(x_{k})\geq0 $ for any $ i $. $\blacksquare$

\end{document}